\documentclass[12pt]{amsart}

\title{Group class operations and homological conditions}
\author{Ioannis Emmanouil and Wei Ren}

\newtheorem{Lemma}{Lemma}[section]
\newtheorem{Proposition}[Lemma]{Proposition}
\newtheorem{Theorem}[Lemma]{Theorem}
\newtheorem{Corollary}[Lemma]{Corollary}

\begin{document}

\begin{abstract}
Kropholler's operation ${\scriptstyle{{\bf LH}}}$
and Talelli's operation $\Phi$ can be often used
to formally enlarge the class of available examples
of groups that satisfy certain homological conditions.
In this paper, we employ this enlargement technique
regarding two specific homological conditions. We
thereby demonstrate the abundance of groups that (a)
have virtually Gorenstein group algebras, as defined
by Beligiannis and Reiten, and (b) satisfy Moore's
conjecture on the relation between projectivity and
relative projectivity, that was studied by Chouinard,
Aljadeff, Cornick, Ginosar, Kropholler and Meir.
\end{abstract}

\maketitle
\tableofcontents

\addtocounter{section}{-1}
\section{Introduction}

\noindent
The study of group actions on spaces is very useful in
cohomological group theory. It offers a geometric intuition
that often clarifies purely algebraic results and also leads
to important constructions that help us understand fine
aspects of the group structure. In order to examine certain
homological finiteness conditions on groups and explore
their relation to complete cohomology, Kropholler has introduced
in \cite{Kro} the ${\scriptstyle{{\bf LH}}}\mathfrak{F}$-groups,
building upon earlier work by Ikenaga \cite{I}. These groups are
constructed hierarchically, starting with the class of finite
groups, and constitute a very big class; it does require some
effort to find groups which are not contained therein. It turns
out that the hierarchical definition of these groups provides
an effective tool for extending homological properties of
finite groups to a much wider class. From the extensive
literature on this topic, we only mention the early
prototypical works \cite{Kro}, \cite{Ben} and \cite{CK2}. More
generally, we may apply Kropholler's construction for any group
class $\mathfrak{C}$ and obtain a group class
${\scriptstyle{{\bf LH}}}\mathfrak{C}$ (which is, in principle,
much bigger that $\mathfrak{C}$). Another example of a group
class operation is provided by Talelli's construction of groups
of type $\Phi$; these groups were defined in \cite{T}, in order
to construct groups that admit a finite dimensional model for the
classifying space for proper actions. Even though the motivation
originates from geometric considerations, the definition of groups
of type $\Phi$ is purely algebraic: It is based on the need to
describe modules of finite projective dimension over the group,
in terms of the finiteness of the projective dimension of their
restriction to finite subgroups. More generally, starting with
any class of groups $\mathfrak{C}$, Talelli's construction provides
us with a group class $\Phi\mathfrak{C}$; see also \cite{MS}.
We may replace the projective dimension with either the flat
or the injective dimension and obtain the analogous group class
operations $\Phi_{flat}$ and $\Phi_{inj}$. An iteration of
Kropholler's operation ${\scriptstyle{{\bf LH}}}$ and Talelli's
operation $\Phi$ (in some of its three forms) often enlarges the
class of available examples of groups that satisfy appropriate
homological conditions.

Examples of this kind of homological conditions can be found
in the realm of Gorenstein homological algebra, the relative
homological theory which is based upon the classes of Gorenstein
projective, Gorenstein flat and Gorenstein injective modules
introduced in \cite{EJ1} and \cite{EJT}. This theory has been
developed rapidly and found interesting applications in algebraic
geometry, representation theory and cohomological group theory.
Even though important advances have been made recently by
\v{S}aroch and \v{S}t$\!$'$\!$ov\'{i}\v{c}ek \cite{SS}, there
are still many properties of the classical homological algebra
notions that are not universally available for their Gorenstein
homological algebra counterparts. The technique of using the
operations ${\scriptstyle{{\bf LH}}}$ and $\Phi$ was employed in
\cite{EmmT2}, in order to exhibit a very big class of groups that
is convenient to work with, as far as these properties are concerned,
in the sense that these properties are satisfied by the Gorenstein
modules over the corresponding group algebras.

In the present work, we provide two more instances of such an
enlargement of the class of examples of groups that satisfy
certain (more or less restrictive) homological conditions. The
first of these examples deals with the notion of virtually
Gorenstein algebras. This notion was introduced for Artin
algebras by Beligiannis and Reiten in \cite{BR}, in order
to study the left-right symmetry of Gorenstein regularity;
cf.\ \cite[Theorem 8.7]{Bel2}. It provides a common
generalization of Gorenstein regular algebras and algebras of
finite representation type. The implications of the virtually
Gorenstein condition to the homotopy (or the stable) categories
of Gorenstein projective and Gorenstein injective modules are
analysed in \cite{Bel1}; see also \cite{C} and \cite{WE}. The
virtually Gorenstein property was studied for arbitrary rings
in \cite{DLW}; it is proved therein that all weakly Gorenstein
regular rings are virtually Gorenstein. Having fixed the
commutative ring $k$, we may consider the class $\mathfrak{W}$
that consists of those groups $G$ for which $kG$ is weakly
Gorenstein regular. If the class $\mathfrak{W}$ is non-empty,
i.e.\ if $k$ is itself weakly Gorenstein regular, then
$\mathfrak{W}$ contains all finite groups. For any
$\mathfrak{W}$-group $G$, \cite[Theorem A]{DLW} implies that
$kG$ is virtually Gorenstein. We now enlarge $\mathfrak{W}$
and let $\overline{\mathfrak{W}}$ be the smallest group class
which contains $\mathfrak{W}$ and is both
${\scriptstyle{{\bf LH}}}$-closed and $\Phi_{inj}$-closed. Then,
$\overline{\mathfrak{W}}$-groups admit a hierarchical description
\`{a} la Kropholler, starting with $\mathfrak{W}$ and using
(transfinite) induction based on the group class operation
$\mathfrak{C} \mapsto {\scriptstyle{{\bf LH}}}\mathfrak{C} \cup
 \Phi_{inj}\mathfrak{C}$;
this process is analysed in $\S $1.IV. The enlargement of the
class of available examples, the central theme of this paper,
is illustrated for the virtually Gorenstein property by the
following result.

\medskip

\noindent
{\bf Theorem A.}
{\em If $G$ is a $\overline{\mathfrak{W}}$-group, then $kG$
is virtually Gorenstein.}

\medskip

\noindent
This result is proved in $\S 4$, using the description of the
$\mbox{Ext}^1$-orthogonal classes of Gorenstein modules over
group algebras of groups that are obtained by performing
Kropholler's operation
${\scriptstyle{{\bf LH}}}$ (see also \cite{Ke}) or Talelli's
operation $\Phi$ to appropriate group classes. The relation
between these operations and $\mbox{Ext}^1$-orthogonality is
analysed in $\S 3$.

Another example of a homological condition that is well-behaved
with respect to the operations ${\scriptstyle{{\bf LH}}}$ and
$\Phi$ is provided by the relation between projectivity and
relative projectivity for modules over group algebras. As an
attempt to generalize the crucial step in the proof of Serre's
theorem on torsion-free groups of finite virtual cohomological
dimension \cite[Chapter VIII, $\S $3]{Bro}, Moore has
conjectured that a module over the group algebra of a group
$G$ must be necessarily projective, provided that its restriction
to a subgroup $H \subseteq G$ of finite index (that satisfies an
additional necessary group-theoretic condition) is projective.
It follows from Chouinard's work \cite{Chou} that finite groups
satisfy Moore's conjecture. This result was extended to
${\scriptstyle{{\bf H}}}_1\mathfrak{F}$-groups in \cite{ACGK}.
More generally, Aljadeff and Meir \cite{AM} have proved that
Moore's conjecture is satisfied by all
${\scriptstyle{{\bf LH}}}\mathfrak{F}$-groups. Building upon
their work, in the present paper we improve the latter result
as follows. Having fixed the commutative ring $k$, let
$\overline{\mathfrak{F}}$ be the smallest class of groups
which contains the class $\mathfrak{F}$ of finite groups
and is ${\scriptstyle{{\bf LH}}}$-closed, $\Phi$-closed
and $\Phi_{flat}$-closed. As before,
$\overline{\mathfrak{F}}$-groups admit a hierarchical
description starting with the class $\mathfrak{F}$ of
finite groups and using (transfinite) induction based
on the group class operation
$\mathfrak{C} \mapsto {\scriptstyle{{\bf LH}}}\mathfrak{C} \cup
 \Phi\mathfrak{C} \cup \Phi_{flat}\mathfrak{C}$.
The following result is proved in $\S 5$.

\medskip

\noindent
{\bf Theorem B.}
{\em If $G$ is an $\overline{\mathfrak{F}}$-group, then $G$
satisfies Moore's conjecture over $k$.}

\medskip

\noindent
Here is a brief description of the contents of the paper:
Following the preliminary Section 1, in Section 2, we fix
the convenient universe of groups within which all Gorenstein
homological algebra notions behave smoothly. In Section 3, we
build on Kendall's work \cite{Ke} and describe the orthogonal
modules with respect to the $\mbox{Ext}^1$-pairing of the
Gorenstein classes over the group algebras of groups that are
obtained by performing Kropholler's or Talelli's operation to
group classes within the universe of groups described in Section
2. As an application of this description, we prove Theorem A
in Section 4. Finally, in Section 5, we discuss Moore's conjecture
and combine the results by Aljadeff and Meir \cite{AM} with certain
simple observations regarding the relation between projectivity
and relative projectivity to prove Theorem B.

\medskip

\noindent
{\em Notations and terminology.}
Let $k$ be a commutative ring and $G$ a group. We shall
denote by ${\tt Proj}(kG)$, ${\tt Flat}(kG)$ and
${\tt Inj}(kG)$ the classes of projective, flat and
injective $kG$-modules respectively. We shall also
denote by $\mathcal{P}(kG)$, $\mathcal{F}(kG)$ and
$\mathcal{I}(kG)$ the classes of $kG$-modules that
appear as kernels of acyclic complexes of projective,
flat or injective $kG$-modules respectively. Finally,
if $M$ is a $kG$-module, then its character (Pontryagin
dual) module is the $kG$-module
$DM = \mbox{Hom}_{\mathbb{Z}}(M,\mathbb{Q}/\mathbb{Z})$
of all additive maps from $M$ to $\mathbb{Q}/\mathbb{Z}$.

\section{Preliminaries}

\noindent
\noindent
In this preliminary section, we record certain basic notions
that are used in the sequel. These notions concern restriction
and (co-)induction of modules over group algebras, basic
properties of Gorenstein modules, cotorsion pairs in module
categories and the hierarchical construction of the closure
of group classes under certain operations.

\vspace{0.1in}

\noindent
{\sc I.\ Restriction and (co-)induction.}
Let $k$ be a commutative ring, $G$ a group and $H \subseteq G$
a subgroup. Restricting the action of $G$ on a $kG$-module $M$
to the subgroup $H$, we obtain the $kH$-module $\mbox{res}_H^GM$.
In this way, we define the restriction functor $\mbox{res}_H^G$
from the category of $kG$-modules to that of $kH$-modules. This
functor has both a left and a right adjoint. Its left adjoint is
the induction functor $\mbox{ind}_H^G$; for any $kH$-module $N$
the induced $kG$-module $\mbox{ind}_H^GN$ is equal to
$kG \otimes_{kH}N$, with $G$ acting on the left of $kG$. The
right adjoint of restriction is the coinduction functor
$\mbox{coind}_H^G$. Here, for any $kH$-module $N$ the coinduced
$kG$-module $\mbox{coind}_H^GN$ is equal to
$\mbox{Hom}_{kH}(kG,N)$, with $G$ acting on the right of $kG$.
Any $kH$-module $N$ is contained as a direct summand in both
$kH$-modules $\mbox{res}_H^G\mbox{ind}_H^GN$ and
$\mbox{res}_H^G\mbox{coind}_H^GN$. For more details, the reader
is referred to \cite[Chapter III, $\S $3]{Bro}. Deriving the two
adjunction isomorphisms, we obtain the Echmann-Shapiro isomorphisms:
If $M$ is a $kG$-module and $N$ a $kH$-module, there are natural
identifications
\[ \mbox{Ext}_{kH}^n \! \left( N,\mbox{res}_H^GM \right) \!
   \simeq
   \mbox{Ext}_{kG}^n \! \left( \mbox{ind}_H^GN,M \right)
   \;\; \mbox{and} \;\;\;
   \mbox{Ext}_{kH}^n \! \left( \mbox{res}_H^GM,N \right) \!
   \simeq
   \mbox{Ext}_{kG}^n \! \left( M,\mbox{coind}_H^GN \right) \]
for all $n \geq 0$. We also note that, in the special case where
$H$ has finite index in $G$, the functors $\mbox{ind}_H^G$ and
$\mbox{coind}_H^G$ are naturally isomorphic; cf.\
\cite[Chapter III, Proposition 5.9]{Bro}.

\vspace{0.1in}

\noindent
{\sc II.\ Gorenstein modules.}
An acyclic complex of projective modules over a ring $R$ is
called totally acyclic if it remains acyclic after applying the
functor $\mbox{Hom}_R(\_\!\_,P)$ for any projective module $P$.
A module $M$ is Gorenstein projective if it is a cokernel of a
totally acyclic complex of projective modules; we denote by
${\tt GProj}(R)$ the class of these modules. Dually, an acyclic
complex of injective modules is called totally acyclic if it
remains acyclic after applying the functor $\mbox{Hom}_R(I,\_\!\_)$
for any injective module $I$. A module $N$ is Gorenstein injective
if it is a kernel of a totally acyclic complex of injective modules;
we denote by ${\tt GInj}(R)$ the class of these modules. An acyclic
complex of flat modules is called totally acyclic if it remains
acyclic after applying the functor $J \otimes_R\_\!\_$ for any
injective right module $J$. A module $L$ is Gorenstein flat if
it is a cokernel of a totally acyclic complex of flat modules;
we denote by ${\tt GFlat}(R)$ the class of these modules. The
Gorenstein modules defined above were introduced by Enochs, Jenda
and Torecillas in \cite{EJ1} and \cite{EJT}; see also Holm's paper
\cite{Hol}. A distinguished class of Gorenstein flat modules has
been introduced by \v{S}aroch and \v{S}t$\!$'$\!$ov\'{i}\v{c}ek
in \cite{SS}: A module $K$ is projectively coresolved Gorenstein
flat if it is a cokernel of an acyclic complex of projective
modules, which remains acyclic after applying the functor
$J \otimes_R\_\!\_$ for any injective right module $J$; the
class of these modules is denoted by ${\tt PGF}(R)$. As shown
in \cite[Theorem 4.4]{SS}, projectively coresolved Gorenstein
flat modules are Gorenstein projective, i.e.\
${\tt PGF}(R) \subseteq {\tt GProj}(R)$.

The Gorenstein projective dimension $\mbox{Gpd}_RM$ of a module
$M$ is the length of a shortest resolution of $M$ by Gorenstein
projective modules. Analogously, we may define the (projectively
coresolved) Gorenstein flat dimension of modules. Dually, the
Gorenstein injective dimension $\mbox{Gid}_RN$ of a module $N$
is the length of a shortest coresolution of $N$ by Gorenstein
injective modules. As shown by Beligiannis \cite[$\S $6]{Bel1},
the finiteness of the Gorenstein projective dimension $\mbox{Gpd}_RM$
of a module $M$ is equivalent to the existence of a complete
projective resolution of $M$, i.e.\ to the existence of a totally
acyclic complex of projective modules that coincides in sufficiently
large degrees with an ordinary projective resolution of $M$. A
characterization of those rings $R$ over which all modules admit
complete projective resolutions is obtained in \cite{CK1}, in
terms of the finiteness of the invariants $\mbox{spli} \, R$ and
$\mbox{silp} \, R$, that were introduced by Gedrich and Gruenberg
in \cite{GG}. Here, $\mbox{spli} \, R$ is the supremum of the
projective dimensions of injective modules and $\mbox{silp} \, R$
is the supremum of the injective dimensions of projective modules.
That result by Cornick and Kropholler was alternatively proved by
Bennis and Mahbou in \cite{BM}, where the notion of the Gorenstein
global dimension of the ring was introduced. More precisely, the
(left) Gorenstein global dimension $\mbox{Ggl.dim} \, R$ of the
ring $R$ is defined as the supremum of the Gorenstein projective
dimension of all modules or, equivalently, as the supremum of the
Gorenstein injective dimension of all modules. Its finiteness is
equivalent to the finiteness of both $\mbox{spli} \, R$ and
$\mbox{silp} \, R$, in which case we have equalities
$\mbox{Ggl.dim} \, R =  \mbox{spli} \, R = \mbox{silp} \, R$.
Following Beligiannis \cite{Bel1}, we say that $R$ is (left)
Gorenstein regular if $\mbox{Ggl.dim} \, R < \infty$.

The description of the finiteness of the Gorenstein weak global
dimension $\mbox{Gwgl.dim} \, R$ of $R$, which is defined as the
supremum of the Gorenstein flat dimension of all modules, is more
difficult to obtain. The homological invariants that play a role
here are $\mbox{sfli} \, R$, the supremum of the flat dimensions
of injective modules, and $\mbox{sfli} \, R^{op}$, the corresponding
invariant involving right modules. As shown by Christensen, Estrada
and Thompson in \cite{CET}, $\mbox{Gwgl.dim} \, R < \infty$ if and
only if both $\mbox{sfli} \, R$ and $\mbox{sfli} \, R^{op}$ are
finite, in which case
$\mbox{Gwgl.dim} \, R =  \mbox{sfli} \, R = \mbox{sfli} \, R^{op}$.
If the ring $R$ is isomorphic with its
opposite $R^{op}$ (for example, if $R=kG$ is the group algebra
of a group $G$ with coefficient in a commutative ring $k$), the
latter result admits a simpler proof; cf.\ \cite[Theorem 5.3]{Emm}.
We say that the ring $R$ is weakly Gorenstein regular if
$\mbox{Gwgl.dim} \, R < \infty$. Any (left or right) Gorenstein
regular ring is weakly Gorenstein regular; in fact, we always
have $\mbox{Gwgl.dim} \, R \leq \mbox{Ggl.dim} \, R$; cf.\
\cite[Theorem 3.7]{WYSZ} or \cite[Theorem 5.13]{CELTW}.

\vspace{0.1in}

\noindent
{\sc III.\ Cotorsion pairs.}
If $R$ is a ring, then the $\mbox{Ext}^1$-pairing induces an
orthogonality relation between modules. In this way, if {\tt S}
is a class of modules, the left orthogonal $^{\perp}{\tt S}$ of
${\tt S}$ is the class consisting of those modules $M$, for which
$\mbox{Ext}^1_R(M,S) = 0$ for all $S \in {\tt S}$. Analogously,
the right orthogonal ${\tt S}^{\perp}$ of {\tt S} is the class
consisting of those modules $N$, for which $\mbox{Ext}^1_R(S,N)=0$
for all $S \in {\tt S}$. If ${\tt U},{\tt V}$ are two module
classes, we say that $({\tt U},{\tt V})$ is a cotorsion pair
(cf.\ \cite{EJ2}) if ${\tt U} = \! \, ^{\perp} {\tt V}$ and
${\tt U}^{\perp} = {\tt V}$. The cotorsion pair is hereditary if
$\mbox{Ext}^i_R(U,V)=0$ for all $i>0$ and all $U \in {\tt U}$ and
$V \in {\tt V}$; it is complete if for any module $M$ there exist
short exact sequences
\[ 0 \longrightarrow V \longrightarrow U \longrightarrow M
     \longrightarrow 0
   \;\;\; \mbox{and } \;\;\;
   0 \longrightarrow M \longrightarrow V' \longrightarrow U'
     \longrightarrow 0 , \]
where $U,U' \in {\tt U}$ and $V,V' \in {\tt V}$. We say that
the cotorsion pair $({\tt U},{\tt V})$ is cogenerated by a class
of modules ${\tt S} \subseteq {\tt U}$ if
${\tt S}^{\perp} = {\tt U}^{\perp}$, so that
${\tt V} = {\tt S}^{\perp}$ and
${\tt U} = {^{\perp}} \! \left( {\tt S}^{\perp} \right)$. As
shown in \cite{ET}, any cotorsion pair which is cogenerated
by a {\em set} of modules is necessarily complete.

Non-trivial examples of cotorsion pairs are provided by
Gorenstein modules. As shown in \cite{SS}, there are complete
hereditary cotorsion pairs
$\left( {\tt PGF}(R),{\tt PGF}(R)^{\perp} \right)$,
$\left( {\tt GFlat}(R),{\tt GFlat}(R)^{\perp} \right)$ and
$\left( ^{\perp}{\tt GInj}(R),{\tt GInj}(R) \right)$ for any
ring $R$. In fact, all three cotorsion pairs are cogenerated
by suitable sets of modules. The situation regarding Gorenstein
projective modules is more subtle: It is shown in \cite{CIS}
that $\left( {\tt GProj}(R),{\tt GProj}(R)^{\perp} \right)$
is a hereditary cotorsion pair as well, but its completeness
appears to depend upon set-theoretic assumptions.

\vspace{0.1in}

\noindent
{\sc IV.\ Continuous operations and the closure of group classes.}
An operation $\mathcal{O}$ on group classes assigns a group class
$\mathcal{O}\mathfrak{C}$ to each group class $\mathfrak{C}$, in
such a way that the following two properties are satisfied:
\newline
(i) $\mathfrak{C} \subseteq \mathcal{O}\mathfrak{C}$ for any group
class $\mathfrak{C}$ and
\newline
(ii) $\mathcal{O}\mathfrak{C}_1 \subseteq \mathcal{O}\mathfrak{C}_2$
for any two group classes $\mathfrak{C}_1,\mathfrak{C}_2$ with
$\mathfrak{C}_1 \subseteq \mathfrak{C}_2$.
\newline
We say that the operation $\mathcal{O}$ is continuous if it also
has the following property:
\newline
(iii) Assume that $\mathfrak{C}_{\alpha}$ is a group class defined
for each ordinal number $\alpha$, such that
$\mathfrak{C}_{\alpha} \subseteq \mathfrak{C}_{\beta}$ for all
$\alpha \leq \beta$, and let $\mathfrak{C}$ be the class consisting
of those groups which are contained in $\mathfrak{C}_{\alpha}$ for
some ordinal number $\alpha$. Then, the group class
$\mathcal{O}\mathfrak{C}$ is the class consisting of those groups
which are contained in $\mathcal{O}\mathfrak{C}_{\alpha}$ for some
$\alpha$.

A group class $\mathfrak{C}$ is said to be closed under an operation
$\mathcal{O}$ if $\mathcal{O}\mathfrak{C} = \mathfrak{C}$. The
$\mathcal{O}$-closure of a group class $\mathfrak{C}$ is the smallest
$\mathcal{O}$-closed group class containing $\mathfrak{C}$; it is the
class consisting of those groups $G$ which are contained in $\mathfrak{D}$
for any $\mathcal{O}$-closed class $\mathfrak{D}$ that contains
$\mathfrak{C}$. If the operation $\mathcal{O}$ is continuous, then
there is an alternative hierarchical (bottom-to-top) description of
the $\mathcal{O}$-closure of a group class $\mathfrak{C}$, that we
now explicit. For any group class $\mathfrak{C}$, we define the
classes $\mathfrak{C}_{\alpha}$ for any ordinal number $\alpha$,
using transfinite induction, as follows:
\newline
(a) $\mathfrak{C}_0 = \mathfrak{C}$,
\newline
(b) $\mathfrak{C}_{\alpha +1} =
     \mathcal{O}\mathfrak{C}_{\alpha}$
     for any ordinal $\alpha$ and
\newline
(c) $\mathfrak{C}_{\alpha} =
     \bigcup_{\beta < \alpha}\mathfrak{C}_{\beta}$
    for any limit ordinal $\alpha$.
\newline
Let $\overline{\mathfrak{C}}$ be the class consisting of those
groups which are contained in $\mathfrak{C}_{\alpha}$, for some
ordinal $\alpha$. It is easily seen that $\overline{\mathfrak{C}}$
is the smallest $\mathcal{O}$-closed group class that contains
$\mathfrak{C}$; in other words, $\overline{\mathfrak{C}}$ is the
$\mathcal{O}$-closure of $\mathfrak{C}$.

If $\mathcal{O}_1, \mathcal{O}_2, \ldots , \mathcal{O}_n$ are
(continuous) operations on group classes, then we may consider
the  (continuous) operation
$\mathcal{O}  = \mathcal{O}_1 \cup \mathcal{O}_2 \cup \ldots
 \cup \mathcal{O}_n$,
which is defined by letting
$\mathcal{O}\mathfrak{C} = \mathcal{O}_1\mathfrak{C} \cup
 \mathcal{O}_2\mathfrak{C} \cup \ldots \cup
 \mathcal{O}_n\mathfrak{C}$
for any group class $\mathfrak{C}$. The $\mathcal{O}$-closure
of $\mathfrak{C}$ is the smallest group class which contains
$\mathfrak{C}$ and is closed under all of the operations
$\mathcal{O}_1, \mathcal{O}_2, \ldots , \mathcal{O}_n$.

We are interested in the continuous operation
${\scriptstyle{{\bf LH}}}$ on group classes that was introduced
by Kropholler in \cite{Kro}. If $\mathfrak{C}$ is any group class,
then we define for each ordinal $\alpha$ the group class
${\scriptstyle{{\bf H}}}_{\alpha}\mathfrak{C}$, using transfinite
induction, as follows: First of all, we let
${\scriptstyle{{\bf H}}}_0\mathfrak{C} = \mathfrak{C}$. For an
ordinal $\alpha >0$, the class
${\scriptstyle{{\bf H}}}_{\alpha}\mathfrak{C}$ consists of those
groups $G$ which admit a cellular action on a finite dimensional
contractible CW-complex $X$, in such a way that each isotropy
subgroup of the action belongs to
${\scriptstyle{{\bf H}}}_{\beta}\mathfrak{C}$ for some ordinal
$\beta < \alpha$ (that may depend upon the particular cell). If
$G,X$ are as above, $k$ is a commutative ring and $M$ is a $kG$-module,
then the cellular chain complex of $X$ induces two exact sequences
of $kG$-modules
\[ 0 \longrightarrow M_d \longrightarrow \ldots
     \longrightarrow M_1 \longrightarrow M_0
     \longrightarrow M \longrightarrow 0 \]
and
\[ 0 \longrightarrow M \longrightarrow M^0
     \longrightarrow M^1 \longrightarrow \ldots
     \longrightarrow M^d \longrightarrow 0 , \]
where $d$ is the dimension of $X$ and the $M_i$'s (resp.\ the
$M^i$'s) are direct sums (resp.\ direct products) of $kG$-modules
of the form $\mbox{ind}_H^G\mbox{res}_H^GM$ (resp.\
$\mbox{coind}_H^G\mbox{res}_H^GM$), for some
${\scriptstyle{{\bf H}}}_{\beta}\mathfrak{C}$-subgroup $H$ of $G$
with $\beta < \alpha$. We say that a group belongs to
${\scriptstyle{{\bf H}}}\mathfrak{C}$ if it belongs to
${\scriptstyle{{\bf H}}}_{\alpha}\mathfrak{C}$ for some $\alpha$.
Then, the class ${\scriptstyle{{\bf LH}}}\mathfrak{C}$ consists
of those groups $G$, all of whose finitely generated subgroups
$H$ are contained in some
${\scriptstyle{{\bf H}}}\mathfrak{C}$-subgroup $K \subseteq G$.
The continuity of the operation ${\scriptstyle{{\bf LH}}}$ is
proved in \cite[Lemma 1.3]{EmmT2}. If the class $\mathfrak{C}$
is subgroup-closed, then the classes
${\scriptstyle{{\bf H}}}\mathfrak{C}$ and
${\scriptstyle{{\bf LH}}}\mathfrak{C}$ are subgroup-closed as
well; in this case, ${\scriptstyle{{\bf LH}}}\mathfrak{C}$
consists of those groups all of whose finitely generated
subgroups are ${\scriptstyle{{\bf H}}}\mathfrak{C}$-groups.

The concept of groups of type $\Phi$, that were introduced by
Talelli in \cite{T}, corresponds to another continuous operation,
applied to the class $\mathfrak{F}$ of finite groups. To be more
precise, if $k$ is a commutative ring and $\mathfrak{C}$ is a
class of groups, then we define $\Phi\mathfrak{C}$ as the class
consisting of those groups $G$, for which the following two
conditions are equivalent for any $kG$-module $M$:
\newline
(i) $\mbox{pd}_{kG}M < \infty$,
\newline
(ii) there is an integer $n$, such that
$\mbox{pd}_{kH}\mbox{res}_H^GM \leq n$ for any
$\mathfrak{C}$-subgroup $H \subseteq G$.
\newline
See also \cite[Definition 5.6]{Bi}. The continuity of $\Phi$
is proved in \cite[Lemma 1.4]{EmmT2}. We shall also consider
the version of the operation $\Phi$ for the flat (resp.\
injective) dimension of modules and define for any group class
$\mathfrak{C}$ the class $\Phi_{flat} \mathfrak{C}$ (resp.\
$\Phi_{inj}\mathfrak{C}$) accordingly.

\section{A convenient class of groups for studying Gorenstein modules}

\noindent
There are a few conjectural properties of the class of Gorenstein
(projective, flat or injective) modules over group algebras that
are only known to hold in particular cases. An example of such a
property is the stability under restriction to subgroups: If $M$
is a Gorenstein (projective, flat or injective) $kG$-module and
$H \subseteq G$ is a subgroup, then is the restricted $kH$-module
$\mbox{res}_H^GM$ also Gorenstein (projective, flat or injective
respectively)? Analogous questions involve the completeness of the
Gorenstein projective cotorsion pair, the relation between
Gorenstein projective and Gorenstein flat modules and the duality
between Gorenstein flat and Gorenstein injective modules.

We fix a commutative ring $k$ and let $\mathfrak{X} = \mathfrak{X}(k)$
be the class consisting of those groups $G$, for which
${\tt PGF}(kG) = {\tt GProj}(kG) = {\mathcal P}(kG)$. In other words,
$G$ is an $\mathfrak{X}$-group if all cokernels of acyclic complexes
of projective $kG$-modules are projectively coresolved Gorenstein flat.
We also consider the class $\mathfrak{X}' = \mathfrak{X}'(k)$ consisting
of those groups $G$, for which ${\tt GFlat}(kG) = {\mathcal F}(kG)$,
and the class $\mathfrak{Y} = \mathfrak{Y}(k)$ consisting of those
groups $G$, for which ${\tt GInj}(kG) = {\mathcal I}(kG)$. In other
words, $G$ is an $\mathfrak{X}'$-group (resp.\ a $\mathfrak{Y}$-group)
if all cokernels (resp.\ kernels) of acyclic complexes of flat (resp.\
injective) $kG$-modules are necessarily Gorenstein flat (resp.\
Gorenstein injective). These classes have been studied in \cite{EmmT2},
where the following result is proved.

\begin{Theorem}
The classes $\mathfrak{X}$, $\mathfrak{X}'$ and $\mathfrak{Y}$
are subgroup-closed and have the following properties:

(i) $\mathfrak{X} = \mathfrak{X}' \supseteq \mathfrak{Y}$; cf.\
\cite[Corollary 2.3]{EmmT2}.

(ii) $\mathfrak{X}$ is closed under the operations
${\scriptstyle{{\bf LH}}}$, $\Phi$ and $\Phi_{flat}$;
cf.\ \cite[Theorem 3.3]{EmmT2}

(iii) $\mathfrak{Y}$ is closed under the operations
${\scriptstyle{{\bf LH}}}$ and $\Phi_{inj}$; cf.\
\cite[Theorem 4.1]{EmmT2}.
\end{Theorem}

\noindent
For later use, we record the following properties of Gorenstein
projective modules over a group algebra. The first three of
these properties are valid over all groups, whereas the latter
three are valid under additional assumptions:
\newline
(1) If $H \subseteq G$ is a subgroup, then induction maps
${\tt GProj}(kH)$ into ${\tt GProj}(kG)$.
\newline
(2) For any $n \geq 1$ any Gorenstein projective $kG$-module
$M$ is an $n$-th syzygy in a projective resolution of another
Gorenstein projective $kG$-module $M'$.
\newline
(3) The orthogonal class ${\tt GProj}(kG)^{\perp}$ is
thick\footnote{We recall that a class of modules is thick
if it is closed under direct summands, extensions, cokernels
of monomorphisms and kernels of epimorphisms.} and contains
all $kG$-modules of finite projective dimension.
\newline
(4) If $H$ is an $\mathfrak{X}$-subgroup of $G$, then
restriction maps ${\tt GProj}(kG)$ into ${\tt GProj}(kH)$.
\newline
(5) If $G$ is an $\mathfrak{X}$-group, then
$\left( {\tt GProj}(kG),{\tt GProj}(kG)^{\perp} \right)$ is a
hereditary cotorsion pair which is cogenerated by a set; in
particular, it is complete.
\newline
(6) If $G$ is an $\mathfrak{X}$-group, then
${\tt GProj}(kG)^{\perp}$ contains all $kG$-modules of
finite flat dimension.
\newline
{\em Proof.}
Property (1) is precisely \cite[Lemma 2.6(i)]{EmmT1}, whereas
(2) and (3) are immediate consequences of the definition of
Gorenstein projectivity; see also \cite[Theorem 2.20]{Hol}.
Property (4) follows since restriction maps $\mathcal{P}(kG)$
into $\mathcal{P}(kH) = {\tt GProj}(kH)$, whereas (5) and (6)
follow since ${\tt GProj}(kG) = {\tt PGF}(kG)$, in view of
\cite[Theorems 4.4 and 4.9]{SS}. \hfill $\Box$

\medskip

\noindent
For completeness, we also list the analogous properties of
Gorenstein injective modules:
\newline
(1') If $H \subseteq G$ is a subgroup, then coinduction maps
${\tt GInj}(kH)$ into ${\tt GInj}(kG)$.
\newline
(2') For any $n \geq 1$ any Gorenstein injective $kG$-module
$M$ is an $n$-th cosyzygy in an injective resolution of another
Gorenstein injective $kG$-module $M'$.
\newline
(3') The orthogonal class $^{\perp}{\tt GInj}(kG)$ is thick
and contains all $kG$-modules of finite injective dimension.
\newline
(4') If $H$ is a $\mathfrak{Y}$-subgroup of $G$, then
restriction maps ${\tt GInj}(kG)$ into ${\tt GInj}(kH)$.
\newline
(5') The pair
$\left( ^{\perp}{\tt GInj}(kG),{\tt GInj}(kG) \right)$ is a
hereditary cotorsion pair which is cogenerated by a set; in
particular, it is complete.
\newline
{\em Proof.}
Property (1') follows by dualizing the argument in
\cite[Lemma 2.6(i)]{EmmT1}, whereas (2) is an immediate
consequence of the definition of Gorenstein injective modules.
Properties (3') and (5') follow from \cite[Theorem 5.6]{SS} and
\cite[Theorem 2.22]{Hol}. Finally, property (4') follows since
restriction maps $\mathcal{I}(kG)$ into
$\mathcal{I}(kH) = {\tt GInj}(kH)$. \hfill $\Box$

\medskip

\noindent
The following simple (and widely known) observation is the
starting point for the description of the classes which are
orthogonal to Gorenstein modules in the next section.

\begin{Corollary}
If $H \subseteq G$ is any subgroup, then restriction maps:

(i) ${\tt GProj}(kG)^{\perp}$ into ${\tt GProj}(kH)^{\perp}$
and

(ii) $^{\perp}{\tt GInj}(kG)$ into $^{\perp}{\tt GInj}(kH)$.
\end{Corollary}
\vspace{-0.05in}
\noindent
{\em Proof.}
Assertion (i) follows from the induction-restriction
(Eckmann-Shapiro) isomorphism, in view of property (1)
above. Analogously, assertion (ii) follows from the
restriction-coinduction (Eckmann-Shapiro) isomorphism,
in view of property (1') above. \hfill $\Box$

\section{$\mbox{Ext}^1$-orthogonality and group class operations}

\noindent
We consider the class $\mathfrak{X}$ defined in the previous
section; in view of Theorem 2.1(ii), $\mathfrak{X}$ is closed
under the operations ${\scriptstyle{{\bf LH}}}$, $\Phi$ and
$\Phi_{flat}$. We note that assertions (i) and (ii) below are
essentially due to Kendall \cite{Ke}.

\begin{Theorem}
Let $\mathfrak{X}_0$ be a subclass of $\mathfrak{X}$. We consider
a group $G$ and a $kG$-module $M$, which is such that
$\mbox{res}_H^GM \in {\tt GProj}(kH)^{\perp}$ for any
$\mathfrak{X}_0$-subgroup $H \subseteq G$.

(i) We have $\mbox{res}_H^GM \in {\tt GProj}(kH)^{\perp}$ for any
${\scriptstyle{{\bf H}}}\mathfrak{X}_0$-subgroup $H \subseteq G$.

(ii) If $\mathfrak{X}_0$ is subgroup-closed, then
$\mbox{res}_H^GM \in {\tt GProj}(kH)^{\perp}$ for any
${\scriptstyle{{\bf LH}}}\mathfrak{X}_0$-subgroup
$H \subseteq G$.

(iii) We have $\mbox{res}_H^GM \in {\tt GProj}(kH)^{\perp}$ for
any $\Phi\mathfrak{X}_0$-subgroup $H \subseteq G$.

(iv) We have $\mbox{res}_H^GM \in {\tt GProj}(kH)^{\perp}$ for
any $\Phi_{flat}\mathfrak{X}_0$-subgroup $H \subseteq G$.
\end{Theorem}
\vspace{-0.05in}
\noindent
{\em Proof.}
(i) We shall proceed by induction on the ordinal $\alpha$, for
which $H \in {\scriptstyle{{\bf H}}}_{\alpha}\mathfrak{X}_0$.
The base case where $\alpha =0$ is precisely our hypothesis on
$M$. For the inductive step, assume that $\alpha >0$ and the
result is true for all
${\scriptstyle{{\bf H}}}_{\beta}\mathfrak{X}_0$-subgroups
of $G$ and all $\beta < \alpha$. We consider an
${\scriptstyle{{\bf H}}}_{\alpha}\mathfrak{X}_0$-subgroup
$H \subseteq G$ and let $N \in {\tt GProj}(kH)$. Then, as
noted in $\S $1.IV, there exists an exact sequence of $kH$-modules
\[ 0 \longrightarrow N_d \longrightarrow \ldots
     \longrightarrow N_1 \longrightarrow N_0
     \longrightarrow N \longrightarrow 0 , \]
where $d$ is the dimension of the $H$-CW-complex witnessing
that $H \in {\scriptstyle{{\bf H}}}_{\alpha}\mathfrak{X}_0$
and each $N_i$ is a direct sum of $kH$-modules of the form
$\mbox{ind}_F^H\mbox{res}_F^HN$, with $F$ an
${\scriptstyle{{\bf H}}}_{\beta}\mathfrak{X}_0$-subgroup of
$H$ for some $\beta < \alpha$. Since
${\scriptstyle{{\bf H}}}_{\beta}\mathfrak{X}_0 \subseteq
 {\scriptstyle{{\bf H}}}_{\beta}\mathfrak{X} \subseteq
 {\scriptstyle{{\bf H}}}\mathfrak{X} \subseteq
 {\scriptstyle{{\bf LH}}}\mathfrak{X} = \mathfrak{X}$,
we may invoke property (4) in Section 2 and conclude that
restriction maps ${\tt GProj}(kH)$ into ${\tt GProj}(kF)$;
in particular, $\mbox{res}_F^HN \in {\tt GProj}(kF)$. As
our induction hypothesis implies that
$\mbox{res}_F^GM \in {\tt GProj}(kF)^{\perp}$, we have
\[ \mbox{Ext}^n_{kH} \! \left( \mbox{ind}_F^H\mbox{res}_F^HN,
   \mbox{res}_H^GM \right) \! = \mbox{Ext}^n_{kF} \! \left(
   \mbox{res}_F^HN,\mbox{res}_F^GM \right) \! = 0 \]
for such a subgroup $F$ and all $n>0$. It follows that
$\mbox{Ext}_{kH}^n \! \left( N_i,\mbox{res}_H^GM \right) \! =0$
for all $i=0,1, \ldots ,d$ and all $n>0$. Using dimension
shifting, we conclude that
$\mbox{Ext}_{kH}^n \! \left( N,\mbox{res}_H^GM \right) \! =0$
for all $n>d$. Since this holds for any $N \in {\tt GProj}(kH)$,
we can easily conclude that the $kH$-module $\mbox{res}_H^GM$ is
actually contained in ${\tt GProj}(kH)^{\perp}$. Indeed, the
Gorenstein projective $kH$-module $N$ is the $d$-th syzygy in
a projective resolution of another Gorenstein projective
$kH$-module $N'$ (this is property (2) in Section 2) and hence
\[ \mbox{Ext}_{kH}^1 \! \left( N,\mbox{res}_H^GM \right) \! =
   \mbox{Ext}_{kH}^{d+1} \! \left( N',\mbox{res}_H^GM \right)
   \! = 0 . \]
This completes the inductive step of the proof, taking care
of any ${\scriptstyle{{\bf H}}}\mathfrak{X}_0$-subgroup $H$
of $G$.

(ii) We now assume that $\mathfrak{X}_0$ is subgroup-closed
and consider an
${\scriptstyle{{\bf LH}}}\mathfrak{X}_0$-subgroup $H$ of $G$.
We proceed by induction on the cardinality $\kappa$ of $H$.
If $\kappa \leq \aleph_0$, then the
${\scriptstyle{{\bf LH}}}\mathfrak{X}_0$-group $H$ is actually
contained in ${\scriptstyle{{\bf H}}}\mathfrak{X}_0$ and we are
done by (i) above. If $\kappa$ is uncountable, then we may
express $H$ as a continuous ascending union of subgroups
$(H_{\lambda})_{\lambda < \kappa}$, each one having cardinality
$< \kappa$. Since the class ${\scriptstyle{{\bf LH}}}\mathfrak{X}_0$
is subgroup-closed, $H_{\lambda}$ is an
${\scriptstyle{{\bf LH}}}\mathfrak{X}_0$-group and our induction
hypothesis implies that
$\mbox{res}^G_{H_{\lambda}}M \in {\tt GProj}(kH_{\lambda})^{\perp}$
for all $\lambda$. We fix
$N \in {\tt GProj}(kH)$ and note that
$\mbox{res}^H_{H_{\lambda}}N \in {\tt GProj}(kH_{\lambda})$ for
all $\lambda$; this follows since
$H_{\lambda} \in {\scriptstyle{{\bf LH}}}\mathfrak{X}_0
 \subseteq {\scriptstyle{{\bf LH}}}\mathfrak{X} =
 \mathfrak{X}$,
in view of property (4) in Section 2. Letting
$N_{\lambda} =
 \mbox{ind}_{H_{\lambda}}^H \mbox{res}^H_{H_{\lambda}}N$,
we compute
\[ \mbox{Ext}^n_{kH} \! \left( N_{\lambda},\mbox{res}_H^GM
   \right) \! = \mbox{Ext}^n_{kH_{\lambda}} \! \left(
   \mbox{res}^H_{H_{\lambda}}N,\mbox{res}_{H_{\lambda}}^GM
   \right) \! = 0 \]
for all $\lambda$ and $n>0$. The continuous ascending family
of subgroups $(H_{\lambda})_{\lambda < \kappa}$ induces a
continuous direct system of $kH$-modules
$(N_{\lambda})_{\lambda < \kappa}$ with surjective structure
maps, whose colimit is $N$. The short exact sequence of
$kH$-modules
\[ 0 \longrightarrow K_{\lambda} \longrightarrow N_0
     \longrightarrow N_{\lambda} \longrightarrow 0 , \]
where $K_{\lambda}$ is the kernel of the structure map
$N_0 \longrightarrow N_{\lambda}$, is the $\lambda$-th term
of a continuous direct system of short exact sequences. The
colimit of the latter direct system of short exact sequences
is the short exact sequence of $kH$-modules
\begin{equation}
 0 \longrightarrow K \longrightarrow N_0
   \longrightarrow N \longrightarrow 0 .
\end{equation}
Here, $K$ is equal to the continuous ascending union of its
submodules $(K_{\lambda})_{\lambda < \kappa}$. Since
$K_{\lambda +1}/K_{\lambda}$ can be identified with
the kernel of the (surjective) structure map
$N_{\lambda} \longrightarrow N_{\lambda +1}$, we conclude that
$\mbox{Ext}^n_{kH}(K_{\lambda +1}/K_{\lambda},\mbox{res}_H^GM)=0$
for all $\lambda$ and $n>0$. Therefore, Eklof's lemma
\cite[Theorem 7.3.4]{EJ2} implies that
$\mbox{Ext}^n_{kH} \! \left( K,\mbox{res}_H^GM \right) \! = 0$
for all $n>0$. Then, the short exact sequence (1) above implies
that
$\mbox{Ext}^n_{kH} \! \left( N,\mbox{res}_H^GM \right) \! = 0$
for all $n>1$. Since this holds for any $N \in {\tt GProj}(kH)$,
we can easily conclude that the $kH$-module $\mbox{res}_H^GM$ is
actually contained in ${\tt GProj}(kH)^{\perp}$. Indeed, the
Gorenstein projective $kH$-module $N$ is the first syzygy in a
projective resolution of another Gorenstein projective $kH$-module
$N'$ (this is property (2) in Section 2) and hence
\[ \mbox{Ext}_{kH}^1 \! \left( N,\mbox{res}_H^GM \right) \! =
   \mbox{Ext}_{kH}^2 \! \left( N',\mbox{res}_H^GM \right) \! = 0 . \]
This completes the proof of (ii).

(iii) We fix a $\Phi\mathfrak{X}_0$-subgroup
$H \subseteq G$. Since
$\Phi\mathfrak{X}_0 \subseteq \Phi\mathfrak{X}
 = \mathfrak{X}$,
we may invoke the completeness of the cotorsion pair
$\left( {\tt GProj}(kH),{\tt GProj}(kH)^{\perp} \right)$
(property (5) in Section 2) to obtain a short exact sequence
of $kH$-modules
\[ 0 \longrightarrow \mbox{res}_H^GM \longrightarrow L
     \longrightarrow N \longrightarrow 0 , \]
where $L \in {\tt GProj}(kH)^{\perp}$ and $N \in {\tt GProj}(kH)$.
Then, for any $\mathfrak{X}_0$-subgroup $F \subseteq H$ we have
an exact sequence of restricted $kF$-modules
\[ 0 \longrightarrow \mbox{res}_F^GM \longrightarrow
     \mbox{res}_F^HL \longrightarrow \mbox{res}_F^HN
     \longrightarrow 0 . \]
Corollary 2.2(i) implies that
$\mbox{res}_F^HL \in {\tt GProj}(kF)^{\perp}$, whereas
property (4) of Section 2 shows that
$\mbox{res}_F^HN \in {\tt GProj}(kF)$. Since
$\mbox{res}_F^GM \in {\tt GProj}(kF)^{\perp}$, in view of
our assumption on $M$, the thickness of the class
${\tt GProj}(kF)^{\perp}$ implies that  $\mbox{res}_F^HN$
is also contained in ${\tt GProj}(kF)^{\perp}$. Since
${\tt GProj}(kF) \cap {\tt GProj}(kF)^{\perp} = {\tt Proj}(kF)$,
it follows that $\mbox{res}_F^HN$ is a projective $kF$-module.
This is the case for any $\mathfrak{X}_0$-subgroup $F$ of the
$\Phi\mathfrak{X}_0$-group $H$ and hence
$\mbox{pd}_{kH}N < \infty$. It follows that
$N \in {\tt GProj}(kH)^{\perp}$ and hence the thickness of
${\tt GProj}(kH)^{\perp}$ implies that
$\mbox{res}_H^GM \in {\tt GProj}(kH)^{\perp}$, as needed.

(iv) We fix a $\Phi_{flat}\mathfrak{X}_0$-subgroup
$H \subseteq G$ and follow verbatim the argument in (iii)
above. In the very last sentence of that proof, we use the
fact that the class ${\tt GProj}(kH)^{\perp}$ contains all
modules of finite flat dimension; the latter claim follows
from property (6) in Section 2, since
$H \in \Phi_{flat}\mathfrak{X}_0 \subseteq
 \Phi_{flat}\mathfrak{X} = \mathfrak{X}$.
\hfill $\Box$

\begin{Corollary}
Let $\mathfrak{X}_0$ be a subgroup-closed subclass of
$\mathfrak{X}$ and consider a group $G$, which is contained
in
${\scriptstyle{{\bf LH}}}\mathfrak{X}_0 \cup
 \Phi\mathfrak{X}_0 \cup \Phi_{flat}\mathfrak{X}_0$.
Then, the following conditions are equivalent for a
$kG$-module $M$:

(i) $M \in {\tt GProj}(kG)^{\perp}$,

(ii) $\mbox{res}_H^GM \in {\tt GProj}(kH)^{\perp}$ for any
$\mathfrak{X}_0$-subgroup $H \subseteq G$.
\end{Corollary}
\vspace{-0.05in}
\noindent
{\em Proof.}
The implication (i)$\rightarrow$(ii) follows from Corollary
2.2(i), whereas the implication (ii)$\rightarrow$(i) follows
from Theorem 3.1. \hfill $\Box$

\vspace{0.1in}

\noindent
We shall prove in the next Section (Proposition 4.3) that
$\mathfrak{X}$ contains all groups $G$ for which $kG$ is
weakly Gorenstein regular; hence, it contains all groups
$G$ for which $kG$ is Gorenstein regular. Therefore, the
following result generalizes \cite[Theorem 2.28]{Ke}: For
any fixed subgroup-closed subclass $\mathfrak{X}_0$ of
$\mathfrak{X}$ and any group $G$ which is obtained by
iterating the operation
${\scriptstyle{{\bf LH}}} \cup \Phi \cup \Phi_{flat}$ on
that subclass, it provides a cogenerating set for the
Gorenstein projective cotorsion pair of $kG$, that is
controlled by the $\mathfrak{X}_0$-subgroups of $G$.

\begin{Proposition}
Let $\mathfrak{X}_0$ be a subgroup-closed subclass of
$\mathfrak{X}$ and assume that for any group
$H \in \mathfrak{X}_0$ the Gorenstein projective cotorsion pair
$\left( {\tt GProj}(kH),{\tt GProj}(kH)^{\perp} \right)$ is
cogenerated by a set ${\tt S}_H$. We also consider the closure
$\overline{\mathfrak{X}_0}$ of $\mathfrak{X}_0$ under the
operation
${\scriptstyle{{\bf LH}}} \cup \Phi \cup \Phi_{flat}$,
as in $\S $1.IV. Then, for any
$\overline{\mathfrak{X}_0}$-group $G$ the Gorenstein projective
cotorsion pair
$\left( {\tt GProj}(kG),{\tt GProj}(kG)^{\perp} \right)$
is cogenerated by the set
${\tt T}_G = \bigcup \left\{ \mbox{ind}_H^G{\tt S}_H :
 H \subseteq G, H \in \mathfrak{X}_0 \right\}$.
\end{Proposition}
\vspace{-0.05in}
\noindent
{\em Proof.}
Groups contained in $\overline{\mathfrak{X}_0}$ admit a
hierarchical description, as explained in $\S $1.IV. More
precisely, we define the classes $\mathfrak{X}_{0,\alpha}$
for any ordinal number $\alpha$, as follows:
\newline
(a) $\mathfrak{X}_{0,0} = \mathfrak{X}_0$,
\newline
(b) $\mathfrak{X}_{0,\alpha +1} =
     {\scriptstyle{{\bf LH}}}\mathfrak{X}_{0,\alpha}
     \cup \Phi\mathfrak{X}_{0,\alpha} \cup
     \Phi_{flat}\mathfrak{X}_{0,\alpha}$
    for any ordinal $\alpha$ and
\newline
(c) $\mathfrak{X}_{0,\alpha} =
     \bigcup_{\beta < \alpha}\mathfrak{X}_{0,\beta}$
    for any limit ordinal $\alpha$.
\newline
Then, $\overline{\mathfrak{X}_0}$ is the class consisting of
those groups $G$, for which $G \in \mathfrak{X}_{0,\alpha}$
for some ordinal $\alpha$. We may therefore prove the claim
in the statement of the Proposition using induction on the least
ordinal $\alpha$, for which  $G \in \mathfrak{X}_{0.\alpha}$.
In view of Theorem 2.1(ii), the class $\overline{\mathfrak{X}_0}$
is a subclass of $\mathfrak{X}$; in particular, this is also the
case for the class $\mathfrak{X}_{0,\alpha}$ for all $\alpha$.

First of all, we assume that
$G \in \mathfrak{X}_{0,0} = \mathfrak{X}_0$. Since
$\mbox{ind}_H^G{\tt S}_H \subseteq
 \mbox{ind}_H^G{\tt GProj}(kH) \subseteq
 {\tt GProj}(kG)$
for any ($\mathfrak{X}_0$-)subgroup $H \subseteq G$, we
have ${\tt T}_G \subseteq {\tt GProj}(kG)$. As $G$ is an
$\mathfrak{X}_0$-subgroup of itself, it is clear that
${\tt S}_G \subseteq {\tt T}_G$. The inclusions
${\tt S}_G \subseteq {\tt T}_G \subseteq {\tt GProj}(kG)$
imply that
${\tt GProj}(kG)^{\perp} \subseteq {\tt T}_G^{\perp}
 \subseteq {\tt S}_G^{\perp}$.
Since ${\tt GProj}(kG)^{\perp} = {\tt S}_G^{\perp}$,
we conclude that
${\tt GProj}(kG)^{\perp} = {\tt T}_G^{\perp}$, as needed.

For the inductive step of the proof, we assume that $G$
is an $\mathfrak{X}_{0,\alpha +1}$-group and the result
is known for all $\mathfrak{X}_{0,\alpha}$-groups. For any
$\mathfrak{X}_{0,\alpha}$-group $K$ our inductive assumption
implies that ${\tt T}_K$ cogenerates the Gorenstein projective
cotorsion pair over $kK$, i.e.\ we have
${\tt T}_K^{\perp} = {\tt GProj}(kK)^{\perp}$. Since
$\mathfrak{X}_{0,\alpha}$ is a subgroup-closed subclass
of $\mathfrak{X}$ and
\[ G \in \mathfrak{X}_{0,\alpha +1} =
   {\scriptstyle{{\bf LH}}}\mathfrak{X}_{0,\alpha}
   \cup \Phi\mathfrak{X}_{0,\alpha} \cup
   \Phi_{flat}\mathfrak{X}_{0,\alpha} , \]
we may apply Corollary 3.2 and conclude that the right orthogonal
class ${\tt GProj}(kG)^{\perp}$ consists of those $kG$-modules
$M$, for which
$\mbox{res}_K^GM \in {\tt GProj}(kK)^{\perp} =
 {\tt T}_K^{\perp}$
for any $\mathfrak{X}_{0,\alpha}$-subgroup $K \subseteq G$.
In view of the induction-restriction (Eckmann-Shapiro) isomorphism,
we have
\[ \mbox{res}_K^GM \in {\tt T}_K^{\perp} \Longleftrightarrow
   M \in \left( \mbox{ind}_K^G{\tt T}_K \right) \! ^{\perp} . \]
It follows that
${\tt GProj}(kG)^{\perp} = \bigcap \left\{ \left(
 \mbox{ind}_K^G{\tt T}_K \right) \! ^{\perp} : K \subseteq G,
 K \in \mathfrak{X}_{0,\alpha} \right\}$
and hence the Gorenstein projective cotorsion pair over $kG$
is cogenerated by
$\, \bigcup \left\{ \mbox{ind}_K^G{\tt T}_K : K \subseteq G,
   K \in \mathfrak{X}_{0,\alpha} \right\}$.
On the other hand, for any $\mathfrak{X}_{0,\alpha}$-subgroup
$K \subseteq G$, we have
\[ \mbox{ind}_K^G{\tt T}_K = {\textstyle{\bigcup}}
   \left\{ \mbox{ind}_K^G\mbox{ind}_H^K{\tt S}_H :
   H \subseteq K , H \in \mathfrak{X}_0 \right\}
   \subseteq {\textstyle{\bigcup}} \left\{
   \mbox{ind}_H^G{\tt S}_H : H \subseteq G,
   H \in \mathfrak{X}_0 \right\} . \]
It follows that
\[ {\textstyle{\bigcup}} \left\{ \mbox{ind}_K^G{\tt T}_K :
   K \subseteq G, K \in \mathfrak{X}_{0,\alpha} \right\}
   \subseteq \, {\textstyle{\bigcup}}
   \left\{ \mbox{ind}_H^G{\tt S}_H : H \subseteq G,
   H \in \mathfrak{X}_0 \right\} = {\tt T}_G . \]
The latter inclusion is actually an equality, since
$\mathfrak{X}_0 \subseteq \mathfrak{X}_{0,\alpha}$
and ${\tt S}_H \subseteq {\tt T}_H$ for any
$\mathfrak{X}_0$-subgroup $H \subseteq G$. Hence,
the Gorenstein projective cotorsion pair over $kG$
is cogenerated by ${\tt T}_G$. \hfill $\Box$

\vspace{0.1in}

\noindent
There are results analogous to Theorem 3.1 and Corollary 3.2,
regarding the left orthogonal to Gorenstein injective modules;
these results will be used in the next section. For completeness,
we shall provide the necessary arguments here. Consider the
class $\mathfrak{Y}$ defined in Section 2; in view of Theorem
2.1(iii), $\mathfrak{Y}$ is closed under the operations
${\scriptstyle{{\bf LH}}}$ and $\Phi_{inj}$.

\begin{Theorem}
Let $\mathfrak{Y}_0$ be a subclass of $\mathfrak{Y}$. We
consider a group $G$ and a $kG$-module $M$, such that
$\mbox{res}_H^GM \in {^{\perp}}{\tt GInj}(kH)$ for any
$\mathfrak{Y}_0$-subgroup $H \subseteq G$.

(i) We have $\mbox{res}_H^GM \in {^{\perp}}{\tt GInj}(kH)$
for any ${\scriptstyle{{\bf H}}}\mathfrak{Y}_0$-subgroup
$H \subseteq G$.

(ii) If $\mathfrak{Y}_0$ is subgroup-closed, then
$\mbox{res}_H^GM \in {^{\perp}}{\tt GInj}(kH)$ for any
${\scriptstyle{{\bf LH}}}\mathfrak{Y}_0$-subgroup
$H \subseteq G$.

(iii) We have $\mbox{res}_H^GM \in {^{\perp}}{\tt GInj}(kH)$
for any $\Phi_{inj}\mathfrak{Y}_0$-subgroup $H \subseteq G$.
\end{Theorem}
\vspace{-0.05in}
\noindent
{\em Proof.}
(i) We use induction on the ordinal $\alpha$, for which
$H \in {\scriptstyle{{\bf H}}}_{\alpha}\mathfrak{Y}_0$. The
base case where $\alpha =0$ is precisely our hypothesis on
$M$. For the inductive step, assume that $\alpha >0$ and the
result is true for all
${\scriptstyle{{\bf H}}}_{\beta}\mathfrak{Y}_0$-subgroups
of $G$ and all $\beta < \alpha$. We consider an
${\scriptstyle{{\bf H}}}_{\alpha}\mathfrak{Y}_0$-subgroup
$H \subseteq G$ and let $N \in {\tt GInj}(kH)$. Then, as
noted in $\S $1.IV, there exists an exact sequence of
$kH$-modules
\[ 0 \longrightarrow N \longrightarrow N^0
     \longrightarrow N^1 \longrightarrow \ldots
     \longrightarrow N^d \longrightarrow 0 , \]
where $d$ is the dimension of the $H$-CW-complex witnessing
that $H \in {\scriptstyle{{\bf H}}}_{\alpha}\mathfrak{Y}_0$
and each $N^i$ is a direct product of $kH$-modules of the
form $\mbox{coind}_F^H\mbox{res}_F^HN$, with $F$ an
${\scriptstyle{{\bf H}}}_{\beta}\mathfrak{Y}_0$-subgroup
of $H$ for some $\beta < \alpha$. Since
${\scriptstyle{{\bf H}}}_{\beta}\mathfrak{Y}_0 \subseteq
 {\scriptstyle{{\bf H}}}_{\beta}\mathfrak{Y} \subseteq
 {\scriptstyle{{\bf H}}}\mathfrak{Y} \subseteq
 {\scriptstyle{{\bf LH}}}\mathfrak{Y} = \mathfrak{Y}$, we
may invoke property (4') in Section 2 and conclude that
restriction maps ${\tt GInj}(kH)$ into ${\tt GInj}(kF)$;
in particular, it follows that
$\mbox{res}_F^HN \in {\tt GInj}(kF)$. Our induction hypothesis
implies that $\mbox{res}_F^GM \in {^{\perp}}{\tt GInj}(kF)$
and hence
\[ \mbox{Ext}^n_{kH} \! \left( \mbox{res}_H^GM,
   \mbox{coind}_F^H\mbox{res}_F^HN \right) =
   \mbox{Ext}^n_{kF} \! \left( \mbox{res}_F^GM,
   \mbox{res}_F^HN \right) = 0 \]
for such a subgroup $F$ and all $n>0$. It follows that
$\mbox{Ext}_{kH}^n \! \left( \mbox{res}_H^GM,N^i \right) =0$
for all $i$ and all $n>0$. Using dimension shifting, we conclude
that $\mbox{Ext}_{kH}^n \! \left( \mbox{res}_H^GM,N \right) =0$
if $n>d$. Since this holds for any $N \in {\tt GInj}(kH)$, we
can easily conclude that the kH-module $\mbox{res}_H^GM$ is
actually contained in $^{\perp}{\tt GInj}(kH)$. Indeed, the
Gorenstein injective $kH$-module $N$ is the $d$-th cosyzygy
in an injective resolution of another Gorenstein injective
$kH$-module $N'$ (this is property (2') in Section 2) and hence
\[ \mbox{Ext}_{kH}^1 \! \left( \mbox{res}_H^GM,N \right) =
   \mbox{Ext}_{kH}^{d+1} \! \left( \mbox{res}_H^GM,N' \right) = 0 . \]
This completes the inductive step of the proof, taking care
of any ${\scriptstyle{{\bf H}}}\mathfrak{Y}_0$-subgroup $H$
of $G$.

(ii) We shall now assume that
$\mathfrak{Y}_0$ is subgroup-closed and consider an
${\scriptstyle{{\bf LH}}}\mathfrak{Y}_0$-subgroup $H$ of $G$.
We proceed by induction on the cardinality $\kappa$ of $H$.
If $\kappa$ is countable, then the
${\scriptstyle{{\bf LH}}}\mathfrak{Y}_0$-group $H$ is actually
contained in ${\scriptstyle{{\bf H}}}\mathfrak{Y}_0$ and we are
done by (i) above. If $\kappa$ is uncountable, then we may
express $H$ as a continuous ascending union of subgroups
$(H_{\lambda})_{\lambda < \kappa}$, each one having cardinality
$< \kappa$. Since the class
${\scriptstyle{{\bf LH}}}\mathfrak{Y}_0$ is subgroup-closed,
$H_{\lambda}$ is an ${\scriptstyle{{\bf LH}}}\mathfrak{Y}_0$-group;
our induction hypothesis implies then that
$\mbox{res}^G_{H_{\lambda}}M \in {^{\perp}}{\tt GInj}(kH_{\lambda})$
for all $\lambda$. We now fix $N \in {\tt GInj}(kH)$ and note that
$\mbox{res}^H_{H_{\lambda}}N \in {\tt GInj}(kH_{\lambda})$ for all
$\lambda$; this follows since
$H_{\lambda} \in {\scriptstyle{{\bf LH}}}\mathfrak{Y}_0
 \subseteq {\scriptstyle{{\bf LH}}}\mathfrak{Y} =
 \mathfrak{Y}$,
in view of property (4') in Section 2. Letting
$M_{\lambda} =
 \mbox{ind}_{H_{\lambda}}^H \mbox{res}^G_{H_{\lambda}}M$,
we note that
\begin{equation}
 \mbox{Ext}^n_{kH} \! \left( M_{\lambda},N \right) =
 \mbox{Ext}^n_{kH_{\lambda}} \! \left(
 \mbox{res}^G_{H_{\lambda}}M,\mbox{res}_{H_{\lambda}}^HN
 \right) = 0
\end{equation}
for all $\lambda$ and $n>0$. The ascending family of subgroups
$(H_{\lambda})_{\lambda < \kappa}$ induces a continuous direct
system of $kH$-modules $(M_{\lambda})_{\lambda < \kappa}$ with
surjective structure maps, whose colimit is $\mbox{res}_H^GM$.
The short exact sequence of $kH$-modules
\[ 0 \longrightarrow T_{\lambda} \longrightarrow M_0
     \longrightarrow M_{\lambda} \longrightarrow 0 , \]
where $T_{\lambda}$ is the kernel of the structure map
$M_0 \longrightarrow M_{\lambda}$, is the $\lambda$-th term
of a continuous direct system of short exact sequences. The
colimit of the latter direct system of short exact sequences
is the short exact sequence of $kH$-modules
\begin{equation}
 0 \longrightarrow T \longrightarrow M_0
   \longrightarrow \mbox{res}_H^GM \longrightarrow 0 .
\end{equation}
Here, $T$ is equal to the continuous ascending union of its
submodules $(T_{\lambda})_{\lambda < \kappa}$. Since
$T_{\lambda +1}/T_{\lambda}$ can be identified with
the kernel of the (surjective) structure map
$M_{\lambda} \longrightarrow M_{\lambda +1}$, we conclude
from (2) that
$\mbox{Ext}^n_{kH} \! \left( T_{\lambda +1}/T_{\lambda},N
 \right) =0$
for all $\lambda$ and $n>0$. Then, Eklof's lemma
\cite[Theorem 7.3.4]{EJ2} implies that
$\mbox{Ext}^n_{kH}(T,N) = 0$ for all $n>0$ and the short
exact sequence (3) above implies that
$\mbox{Ext}^n_{kH} \! \left( \mbox{res}_H^GM,N \right) \! =0$
for all $n>1$. Since this holds for any $N \in {\tt GInj}(kH)$,
we can easily conclude that the $kH$-module $\mbox{res}_H^GM$
is actually contained in ${^{\perp}}{\tt GInj}(kH)$. Indeed,
the Gorenstein injective $kH$-module $N$ is the first cosyzygy
in an injective resolution of another Gorenstein injective
$kH$-module $N'$ (this is  property (2') in Section 2) and
hence
\[ \mbox{Ext}_{kH}^1 \! \left( \mbox{res}_H^GM,N \right) \! =
   \mbox{Ext}_{kH}^2 \! \left( \mbox{res}_H^GM,N' \right) \!
   = 0 . \]
This completes the proof of (ii).\footnote{One could replace
the variation of Benson's lemma \cite[Lemma 5.6]{Ben} in the
argument above and use instead the closure of the left
orthogonal class ${^{\perp}}{\tt GInj}(kG)$ under colimits
(cf.\ \cite[Proposition 3.1]{Gil}), in order to conclude from
(2) that the colimit $\mbox{res}_H^GM$ of the direct system
$(M_{\lambda})_{\lambda < \kappa}$ is also contained in
${^{\perp}}{\tt GInj}(kG)$.}

(iii) We fix a $\Phi_{inj}\mathfrak{Y}_0$-subgroup
$H \subseteq G$. Since the cotorsion pair
$\left( {^{\perp}}{\tt GInj}(kH),{\tt GInj}(kH) \right)$ is
complete (property (5') in Section 2), there exists a short
exact sequence of $kH$-modules
\[ 0 \longrightarrow N \longrightarrow L
     \longrightarrow \mbox{res}_H^GM \longrightarrow 0 , \]
where $L \in {^{\perp}}{\tt GInj}(kH)$ and $N \in {\tt GInj}(kH)$.
Then, for any $\mathfrak{Y}_0$-subgroup $F \subseteq H$ we obtain
the exact sequence of restricted $kF$-modules
\[ 0 \longrightarrow \mbox{res}_F^HN \longrightarrow
     \mbox{res}_F^HL \longrightarrow \mbox{res}_F^GM
     \longrightarrow 0 . \]
Corollary 2.2(ii) implies that
$\mbox{res}_F^HL \in {^{\perp}}{\tt GInj}(kF)$, whereas
property (4') of Section 2 shows that
$\mbox{res}_F^HN \in {\tt GInj}(kF)$. Since
$\mbox{res}_F^GM \in {^{\perp}}{\tt GInj}(kF)$, in view of
our assumption on $M$, the thickness of the class
${^{\perp}}{\tt GInj}(kF)$ implies that  $\mbox{res}_F^HN$
is also contained in ${^{\perp}}{\tt GInj}(kF)$. Since
${\tt GInj}(kF) \cap {^{\perp}}{\tt GInj}(kF) = {\tt Inj}(kF)$,
it follows that $\mbox{res}_F^HN$ is an injective $kF$-module.
This is the case for any $\mathfrak{Y}_0$-subgroup $F$ of the
$\Phi_{inj}\mathfrak{Y}_0$-group $H$ and hence
$\mbox{id}_{kH}N < \infty$. It follows that
$N \in {^{\perp}}{\tt GInj}(kH)$ and hence the thickness of
${^{\perp}}{\tt GInj}(kH)$ implies that
$\mbox{res}_H^GM \in {^{\perp}}{\tt GInj}(kH)$, as needed.
\hfill $\Box$

\begin{Corollary}
Let $\mathfrak{Y}_0$ be a subgroup-closed subclass
of $\mathfrak{Y}$ and consider a group $G$, which
is contained in
${\scriptstyle{{\bf LH}}}\mathfrak{Y}_0 \cup
 \Phi_{inj}\mathfrak{Y}_0$.
Then, the following are equivalent for a $kG$-module $M$:

(i) $M \in {^{\perp}}{\tt GInj}(kG)$,

(ii) $\mbox{res}_H^GM \in {^{\perp}}{\tt GInj}(kH)$ for any
$\mathfrak{Y}_0$-subgroup $H \subseteq G$.
\end{Corollary}
\vspace{-0.05in}
\noindent
{\em Proof.}
The implication (i)$\rightarrow$(ii) follows from Corollary
2.2(ii), whereas the implication (ii)$\rightarrow$(i) follows
from Theorem 3.4. \hfill $\Box$

\section{Virtually Gorenstein group algebras}

\noindent
We fix a commutative ring $k$ and consider groups $G$ for
which the group algebra $kG$ is virtually Gorenstein, i.e.\
for which ${\tt GProj}(kG)^{\perp} = {^{\perp}}{\tt GInj}(kG)$.
This homological condition on $kG$ has interesting consequences
regarding the (homotopy) properties of Gorenstein modules.
As an indication of these consequences, we note that if $kG$
is virtually Gorenstein, then:\footnote{Properties (a), (b),
(c) and (d) are actually valid over any (left) virtually
Gorenstein ring $R$.}
\newline
(a) (cf.\ \cite[Proposition 3.16]{WE}) There is an equality
${\tt GProj}(kG) = {\tt PGF}(kG)$; hence,
\cite[Theorem 4.4]{SS} implies that
${\tt GProj}(kG) \subseteq {\tt GFlat}(kG)$ and
\cite[Theorem 4.9]{SS} implies that the cotorsion pair
$\left( {\tt GProj}(kG),{\tt GProj}(kG)^{\perp} \right)$
is complete.
\newline
(b) (cf.\ \cite[Proposition 2.6]{C}) Let $M$ be a $kG$-module.
Any resolution of $M$ which is obtained by iterating special
Gorenstein projective precovers remains acyclic by applying
the functor $\mbox{Hom}_{kG}(\_\!\_,N)$ for any Gorenstein
injective $kG$-module $N$. Dually, any coresolution of $M$
which is obtained by iterating special Gorenstein injective
preenvelopes remains acyclic by applying the functor
$\mbox{Hom}_{kG}(L,\_\!\_)$ for any Gorenstein projective
$kG$-module $L$.
\newline
(c) (cf.\ \cite[Proposition 2.2]{C}) An acyclic complex of
$kG$-modules remains acyclic by applying the functor
$\mbox{Hom}_{kG}(\_\!\_,N)$ for any Gorenstein injective
$kG$-module $N$ if and only if it remains acyclic by
applying the functor $\mbox{Hom}_{kG}(L,\_\!\_)$ for any
Gorenstein projective $kG$-module $L$.
\newline
(d) (cf.\ \cite[Example 3.18(2)]{WE}) Let $\mathcal{C}$
be the class of those acyclic complexes of $kG$-modules
that satisfy the equivalent conditions described in (c)
above. Then, the relative derived category, which is
obtained from the homotopy category of $kG$ by inverting
those chain maps whose mapping cone is a $\mathcal{C}$-complex,
is triangulated equivalent to certain homotopy categories
of Gorenstein projective (resp.\ Gorenstein injective)
$kG$-modules.

\medskip

\noindent
We consider the classes $\mathfrak{X}$ and $\mathfrak{Y}$
defined in Section 2 and let $\mathfrak{Z} = \mathfrak{Z}(k)$
be the class consisting of those $\mathfrak{Y}$-groups $G$,
for which $kG$ is virtually Gorenstein.

\begin{Lemma}
The class $\mathfrak{Z}$ is subgroup-closed.
\end{Lemma}
\vspace{-0.05in}
\noindent
{\em Proof.}
Let $G$ be a $\mathfrak{Z}$-group and consider a subgroup
$H \subseteq G$. Then, $G \in \mathfrak{Y}$ and $kG$ is
virtually Gorenstein. Since $\mathfrak{Y}$ is
subgroup-closed (cf.\ Theorem 2.1), it follows that
$H \in \mathfrak{Y}$. It remains to show that $kH$ is
virtually Gorenstein.

Consider a $kH$-module $M \in {\tt GProj}(kH)^{\perp}$. Since
$H \in \mathfrak{Y} \subseteq \mathfrak{X}$ (cf.\ Theorem 2.1(i)),
restriction maps ${\tt GProj}(kG)$ into ${\tt GProj}(kH)$; this
is precisely property (4) in Section 2. It follows that the
coinduced $kG$-module
$\mbox{coind}_H^GM$ is contained in ${\tt GProj}(kG)^{\perp}$,
i.e.\ in ${^{\perp}}{\tt GInj}(kG)$. Invoking Corollary 2.2(ii),
we may now conclude that
$\mbox{res}_H^G\mbox{coind}_H^GM \in {^{\perp}}{\tt GInj}(kH)$.
Since $M$ is a direct summand of
$\mbox{res}_H^G\mbox{coind}_H^GM$, it follows that
$M \in {^{\perp}}{\tt GInj}(kH)$. Hence, we have shown that
${\tt GProj}(kH)^{\perp} \subseteq {^{\perp}}{\tt GInj}(kH)$.

Conversely, assume that $M \in {^{\perp}}{\tt GInj}(kH)$. Since
$H$ is a $\mathfrak{Y}$-group, restriction maps ${\tt GInj}(kG)$
into ${\tt GInj}(kH)$; this is precisely property (4') in Section
2. It follows that the induced $kG$-module
$\mbox{ind}_H^GM$ is contained in ${^{\perp}}{\tt GInj}(kG)$,
i.e.\ in ${\tt GProj}(kG)^{\perp}$. Corollary 2.2(i) then
implies that
$\mbox{res}_H^G\mbox{ind}_H^GM \in {\tt GProj}(kH)^{\perp}$.
Since $M$ is a direct summand of $\mbox{res}_H^G\mbox{ind}_H^GM$,
it follows that $M \in {\tt GProj}(kH)^{\perp}$. We have thus
shown that
${^{\perp}}{\tt GInj}(kH) \subseteq {\tt GProj}(kH)^{\perp}$.
\hfill $\Box$

\vspace{0.1in}

\noindent
The following result will provide us with an abundance of
examples of virtually Gorenstein group algebras.

\begin{Theorem}
The class $\mathfrak{Z}$ is closed under the operations
(i) ${\scriptstyle{{\bf LH}}}$ and (ii) $\Phi_{inj}$.
\end{Theorem}
\vspace{-0.05in}
\noindent
{\em Proof.}
(i) We have to show that
${\scriptstyle{{\bf LH}}}\mathfrak{Z} = \mathfrak{Z}$, i.e.\
that ${\scriptstyle{{\bf LH}}}\mathfrak{Z} \subseteq \mathfrak{Z}$.
To that end, we consider a group
$G \in {\scriptstyle{{\bf LH}}}\mathfrak{Z}$ and aim at
proving that $G \in \mathfrak{Z}$. We note that
${\scriptstyle{{\bf LH}}}\mathfrak{Z} \subseteq
 {\scriptstyle{{\bf LH}}}\mathfrak{Y} = \mathfrak{Y}$
and hence $G \in \mathfrak{Y}$. It remains to show that $kG$
is virtually Gorenstein. For any $\mathfrak{Z}$-subgroup
$H \subseteq G$ we denote by
$\left( \mbox{res}_H^G \right)^{-1} {\tt GProj}(kH)^{\perp}$
the class of those $kG$-modules $M$ for which
$\mbox{res}_H^GM \in {\tt GProj}(kH)^{\perp}$ and apply
Corollary 3.2 (with
$\mathfrak{X}_0 = \mathfrak{Z} \subseteq \mathfrak{Y}
 \subseteq \mathfrak{X}$
therein), in order to conclude that
\[ {\tt GProj}(kG)^{\perp} =
   {\textstyle{\bigcap}} \left\{
   \left( \mbox{res}_H^G \right)^{-1} {\tt GProj}(kH)^{\perp}
   : H \subseteq G , H \in \mathfrak{Z} \right\} . \]
Analogously, for any $\mathfrak{Z}$-subgroup $H \subseteq G$
we denote by
$\left( \mbox{res}_H^G \right)^{-1} {^{\perp}}{\tt GInj}(kH)$
the class of those $kG$-modules $M$ for which
$\mbox{res}_H^GM \in {^{\perp}}{\tt GInj}(kH)$ and apply
Corollary 3.5 (with
$\mathfrak{Y}_0 = \mathfrak{Z} \subseteq \mathfrak{Y}$ therein),
in order to conclude that
\[ {^{\perp}}{\tt GInj}(kG) =
   {\textstyle{\bigcap}} \left\{
   \left( \mbox{res}_H^G \right)^{-1} {^{\perp}}{\tt GInj}(kH)
   : H \subseteq G , H \in \mathfrak{Z} \right\} . \]
Since ${\tt GProj}(kH)^{\perp} = {^{\perp}}{\tt GInj}(kH)$ for
any $\mathfrak{Z}$-subgroup $H \subseteq G$, it follows that
${\tt GProj}(kG)^{\perp} = {^{\perp}}{\tt GInj}(kG)$. The group
algebra $kG$ is therefore virtually Gorenstein and hence
$G \in \mathfrak{Z}$.

(ii) In order to show that
$\Phi_{inj}\mathfrak{Z} = \mathfrak{Z}$, i.e.\ that
$\Phi_{inj}\mathfrak{Z} \subseteq \mathfrak{Z}$, we consider
a $\Phi_{inj}\mathfrak{Z}$-group $G$ and aim at proving that
$G \in \mathfrak{Z}$. First of all, we note that
$\Phi_{inj}\mathfrak{Z} \subseteq \Phi_{inj}\mathfrak{Y} =
 \mathfrak{Y}$,
so that $G \in \mathfrak{Y}$. Since $G$ is contained in
$\Phi_{inj}\mathfrak{Z} \subseteq \Phi_{flat}\mathfrak{Z}$
(cf.\ \cite[Lemma 1.5]{EmmT2}), we can apply Corollaries 3.2
and 3.5 in exactly the same way as in (i) above and show that
${\tt GProj}(kG)^{\perp} = {^{\perp}}{\tt GInj}(kG)$, i.e.\
that the group algebra $kG$ is virtually Gorenstein. \hfill $\Box$

\vspace{0.1in}

\noindent
Let $\mathfrak{W} = \mathfrak{W}(k)$ be the class consisting
of those groups $G$, for which $kG$ is weakly Gorenstein regular;
in other words, $G \in \mathfrak{W}$ if and only if
$\mbox{sfli} \, kG < \infty$. This class contains all groups $G$,
for which $kG$ is Gorenstein regular. We note that $\mathfrak{W}$
is subgroup-closed; indeed, for any subgroup $H$ of a group $G$
there is an inequality $\mbox{sfli} \, kH \leq \mbox{sfli} \, kG$.
Hence, $\mathfrak{W}$ is non-empty if and only if it contains
the trivial group, i.e.\ if and only if $k$ is itself weakly
Gorenstein regular. In that case, $\mathfrak{W}$ is known to
contain the class $\mathfrak{F}$ of finite groups; see, for
example, \cite{RY}.

\begin{Proposition}
The class $\mathfrak{W}$ is a subclass of $\mathfrak{Z}$.
\end{Proposition}
\vspace{-0.05in}
\noindent
{\em Proof.}
We have to show that any $\mathfrak{W}$-group $G$ has the
following two properties:
(i) $G \in \mathfrak{Y}$, i.e.\ ${\tt GInj}(kG) = \mathcal{I}(kG)$
and
(ii) $kG$ is virtually Gorenstein.

Since $\mbox{sfli} \, kG < \infty$, all injective $kG$-modules
have finite flat dimension. On the other hand, if $I^*$ is any
acyclic complex of injective $kG$-modules and $M$ is any $kG$-module
of finite flat dimension, then the complex of abelian groups
$\mbox{Hom}_{kG}(M,I^*)$ is acyclic; cf.\
\cite[Corollary 5.9]{Sto}. This shows that $G$ satisfies property
(i) above; see also \cite[Theorem 2.6]{WZ}. As far as property (ii)
is concerned, this is precisely \cite[Theorem A]{DLW}; see also
\cite[Lemma 3.2]{WZ}. \hfill $\Box$

\vspace{0.1in}

\noindent
Let $\overline{\mathfrak{W}}$ be the smallest class of groups
which contains $\mathfrak{W}$ and is
${\scriptstyle{{\bf LH}}}$-closed and $\Phi_{inj}$-closed.
Groups in $\overline{\mathfrak{W}}$ admit a hierarchical
description, as explained in $\S 1$.IV. The following Corollary
implies Theorem A in the Introduction; it is itself an immediate
consequence of Theorem 4.2 and Proposition 4.3.

\begin{Corollary}
The class $\overline{\mathfrak{W}}$ is contained in
$\mathfrak{Z}$; in particular, $kG$ is virtually
Gorenstein for any $\overline{\mathfrak{W}}$-group $G$.
\end{Corollary}

\section{Groups satisfying Moore's conjecture}

\noindent
Motivated by Serre's theorem on torsion-free groups of finite
virtual cohomological dimension \cite[Chapter VIII, $\S $3]{Bro},
Moore has conjectured that a module $M$ over the group algebra
of a group $G$ is projective, whenever the restriction of $M$
to a subgroup $H \subseteq G$ of finite index that satisfies
an additional group-theoretic condition is projective. More
precisely, we say that the triple $(G,H,k)$, where $G$ is a
group and $H \subseteq G$ is a subgroup of finite index,
satisfies Moore's condition if for any element
$g \in G \setminus H$ one of the following two conditions
hold:

(i) the intersection $\langle g \rangle \cap H$ is non-trivial
or

(ii) $g$ has finite order, which is invertible in $k$.
\newline
Following \cite{AM}, we say that $G$ satisfies Moore's conjecture
over $k$ if for any subgroup $H \subseteq G$ of finite index,
such that the triple $(G,H,k)$ satisfies Moore's condition,
and any $kG$-module $M$, for which
$\mbox{res}_H^GM \in {\tt Proj}(kH)$, we actually have
$M \in {\tt Proj}(kG)$. We denote by
$\mathfrak{M} = \mathfrak{M}(k)$ the class consisting of those
groups that satisfy Moore's conjecture over $k$. Replacing
projective modules by flat (resp.\ injective) modules in the
above discussion, we obtain versions of Moore's conjecture for
flat (resp.\ injective) modules; cf.\ \cite{BS1}. We shall
denote by $\mathfrak{M}_{flat} = \mathfrak{M}_{flat}(k)$ and
$\mathfrak{M}_{inj} = \mathfrak{M}_{inj}(k)$ the corresponding
classes of groups. In other words, $\mathfrak{M}_{flat}$ (resp.\
$\mathfrak{M}_{inj}$) is the class consisting of those groups
that satisfy the flat (resp.\ injective) version of Moore's
conjecture over $k$.

\begin{Proposition}
There are inclusions
$\mathfrak{M}_{inj} \subseteq \mathfrak{M}_{flat} \subseteq
 \mathfrak{M}$.
\end{Proposition}
\vspace{-0.05in}
\noindent
{\em Proof.}
The first inclusion $\mathfrak{M}_{inj} \subseteq \mathfrak{M}_{flat}$
is a consequence of Lambek's flatness criterion \cite{L}, as noted
in \cite[Proposition 2.6]{BS1}. The inclusion
$\mathfrak{M}_{flat} \subseteq \mathfrak{M}$ follows from the
Benson-Goodearl theorem \cite[Corollary 4.8]{BG}, as in the proof
of \cite[Proposition 2.7]{BS1}. \hfill $\Box$

\vspace{0.1in}

\noindent
The following simple (and well-known) result will be helpful
in the proof of Theorems 5.3-5.5 below. We include a proof
for the sake of completeness.

\begin{Proposition}
Let $G$ be a group, $H \subseteq G$ a subgroup of finite index
and $M$ a $kG$-module.

(i) If $\mbox{pd}_{kG}M < \infty$ and
$\mbox{res}_H^GM \in {\tt Proj}(kH)$, then $M$ is projective.

(ii) If $\mbox{id}_{kG}M < \infty$ and
$\mbox{res}_H^GM \in {\tt Inj}(kH)$, then $M$ is injective.

(iii) If $\mbox{fd}_{kG}M < \infty$ and
$\mbox{res}_H^GM \in {\tt Flat}(kH)$, then $M$ is flat.

(iv) If $\mbox{fd}_{kG}M < \infty$ and
$\mbox{res}_H^GM \in {\tt Proj}(kH)$, then $M$ is projective.
\end{Proposition}
\vspace{-0.05in}
\noindent
{\em Proof.}
(i) Assume that $n = \mbox{pd}_{kG}M \geq 1$ and consider an
exact sequence of $kG$-modules
\[ 0 \longrightarrow P_n
     \stackrel{f}{\longrightarrow} P_{n-1}
     \longrightarrow \ldots \longrightarrow P_0
     \longrightarrow M \longrightarrow 0 , \]
where $P_0 , \ldots , P_{n-1} , P_n \in {\tt Proj}(kG)$.
We denote this exact sequence by $X_*$. Since $f$ is not
a split monomorphism, the complex of abelian groups
$\mbox{Hom}_{kG}(X_*,P_n)$ is not acyclic. We shall reach
a contradiction, by showing that $\mbox{Hom}_{kG}(X_*,P)$
is acyclic for any projective $kG$-module $P$. In fact,
it suffices to consider the case where $P$ is $kG$-free.
Then, $P = \mbox{ind}_H^GQ = \mbox{coind}_H^GQ$ for a
suitable (free) $kH$-module $Q$ and hence
\[ \mbox{Hom}_{kG}(X_*,P) =
   \mbox{Hom}_{kG} \! \left( X_*,\mbox{coind}_H^GQ \right) \! =
   \mbox{Hom}_{kH} \! \left( \mbox{res}_H^GX_*,Q \right) \! . \]
In view of our assumption on $\mbox{res}_H^GM$, the restricted
complex $\mbox{res}_H^GX_*$ is contractible (being a right-bounded
acyclic complex of projective $kH$-modules). Therefore, the complex
of abelian groups
$\mbox{Hom}_{kG}(X_*,P) =
 \mbox{Hom}_{kH} \! \left( \mbox{res}_H^GX_*,Q \right)$
is contractible as well; in particular, it is acyclic.

(ii) Assume that $n = \mbox{id}_{kG}M \geq 1$ and consider an
exact sequence of $kG$-modules
\[ 0 \longrightarrow M \longrightarrow I^0
     \longrightarrow \ldots \longrightarrow I^{n-1}
     \stackrel{g}{\longrightarrow} I^n \longrightarrow 0 , \]
where $I^0 , \ldots , I^{n-1} , I^n \in {\tt Inj}(kG)$.
We denote this exact sequence by $Y^*$. Since $g$ is not
a split epimorphism, the complex of abelian groups
$\mbox{Hom}_{kG}(I^n,Y^*)$ is not acyclic. We shall reach
a contradiction, by showing that $\mbox{Hom}_{kG}(I,Y^*)$
is acyclic for any injective $kG$-module $I$. In fact, it
suffices to consider the case where $I$ is a coinduced
module of the form $\mbox{coind}_1^GA$, for a divisible abelian
group $A$. Then, $I = \mbox{coind}_H^GJ = \mbox{ind}_H^GJ$
for a suitable (injective) $kH$-module $J$ and hence
\[ \mbox{Hom}_{kG}(I,Y^*) =
   \mbox{Hom}_{kG} \! \left( \mbox{ind}_H^GJ,Y^* \right) \! =
   \mbox{Hom}_{kH} \! \left( J,\mbox{res}_H^GY^* \right) \! . \]
In view of our assumption on $\mbox{res}_H^GM$, the restricted
complex $\mbox{res}_H^GY^*$ is contractible (being a left-bounded
acyclic complex of injective $kH$-modules). It follows that the
complex of abelian groups
$\mbox{Hom}_{kG}(I,Y^*) =
 \mbox{Hom}_{kH} \! \left( J,\mbox{res}_H^GY^* \right)$
is contractible as well; in particular, it is acyclic.

(iii) This follows from assertion (ii) above, in view of
Lambek's flatness criterion \cite{L}. Indeed, the Pontryagin
dual $DM$ of $M$ has injective dimension
$\mbox{id}_{kG}DM = \mbox{fd}_{kG}M < \infty$, whereas the
restricted $kH$-module $\mbox{res}_H^GDM = D\mbox{res}_H^GM$
is injective. We may then invoke (ii) and conclude that
$DM \in {\tt Inj}(kG)$, so that $M \in {\tt Flat}(kG)$.

(iv) Assume that $\mbox{fd}_{kG}M=n$ and let $N = \Omega^nM$
be the $n$-th syzygy in a projective resolution of $M$. Then,
the $kG$-module $N$ is flat and its restriction
$\mbox{res}_H^GN = \Omega^n\mbox{res}_H^GM$ is $kH$-projective.
Invoking the Benson-Goodearl theorem \cite[Corollary 4.8]{BG},
we conclude that $N \in {\tt Proj}(kG)$ and hence
$\mbox{pd}_{kG}M \leq n < \infty$. Then, the projectivity
of $M$ is a consequence of (i) above. \hfill $\Box$

\vspace{0.1in}

\noindent
The next results are based upon and, at the same time,
complement previous work on Moore's conjecture by
Aljadeff, Meir, Bahlekeh and Salarian; cf.\ \cite{AM},
\cite{BS1}.

\begin{Theorem}
Let $G$ be a group and $H \subseteq G$ a subgroup of finite
index, such that the triple $(G,H,k)$ satisfies Moore's
condition. If $M$ is a $kG$-module with
$\mbox{res}_H^GM \in {\tt Proj}(kH)$, then:

(i) $\mbox{res}_S^GM \in {\tt Proj}(kS)$ for any
${\scriptstyle{{\bf LH}}}\mathfrak{M}$-subgroup
$S \subseteq G$,

(ii) $\mbox{res}_S^GM \in {\tt Proj}(kS)$ for any
$\Phi\mathfrak{M}$-subgroup $S \subseteq G$ and

(iii) $\mbox{res}_S^GM \in {\tt Proj}(kS)$ for any
$\Phi_{flat}\mathfrak{M}$-subgroup $S \subseteq G$.
\end{Theorem}
\vspace{-0.05in}
\noindent
{\em Proof.}
For the subgroup $S$ under consideration, we let
$N = \mbox{res}_S^GM$ and note that
$\mbox{res}_{H \cap S}^SN = \mbox{res}_{H \cap S}^GM =
 \mbox{res}_{H \cap S}^H\mbox{res}_H^GM$.
Then, being a restriction of the projective $kH$-module
$\mbox{res}_H^GM$, the latter $k(H \cap S)$-module is
projective as well. We also note that $H \cap S$ has
finite index in $S$ and the triple $(S,H \cap S,k)$
satisfies Moore's condition.

(i) We shall first consider the case where $S$ is an
${\scriptstyle{{\bf H}}}\mathfrak{M}$-subgroup of $G$ and
proceed by induction on the ordinal $\alpha$, for which
$S \in {\scriptstyle{{\bf H}}}_{\alpha}\mathfrak{M}$,
following the argument in the proof of \cite[Proposition 4.2]{AM}.
If $\alpha =0$, then $S \in \mathfrak{M}$ and hence $S$
satisfies Moore's conjecture over $k$. The $kS$-module $N$
is therefore projective (as its restriction to $H \cap S$
is projective). For the inductive step, assume that
$\alpha >0$ and the result is true for all
${\scriptstyle{{\bf H}}}_{\beta}\mathfrak{M}$-subgroups
of $G$ and all $\beta < \alpha$. As we noted in $\S $1.IV,
there exists a short exact sequence of $kS$-modules
\[ 0 \longrightarrow N_d \longrightarrow \ldots
     \longrightarrow N_1 \longrightarrow N_0
     \longrightarrow N \longrightarrow 0 , \]
where $d$ is the dimension of the $S$-CW-complex witnessing
that $S \in {\scriptstyle{{\bf H}}}_{\alpha}\mathfrak{M}$
and each $N_i$ is a direct sum of $kS$-modules of the form
$\mbox{ind}_F^S\mbox{res}_F^SN$, with $F$ an
${\scriptstyle{{\bf H}}}_{\beta}\mathfrak{M}$-subgroup of
$S$ for some $\beta < \alpha$. Our inductive assumption
implies that the $kF$-module
$\mbox{res}_F^SN = \mbox{res}_F^GM$ is projective, so that
$\mbox{ind}_F^S\mbox{res}_F^SN \in {\tt Proj}(kS)$
for all such subgroups $F$. Then, the $kS$-module $N_i$ is
projective for all $i=0,1, \ldots ,d$ and hence
$\mbox{pd}_{kS}N \leq d$. Invoking Proposition 5.2(i), we
conclude that the $kS$-module $N$ is actually projective.

We now assume that $S$ is an
${\scriptstyle{{\bf LH}}}\mathfrak{M}$-subgroup of $G$. If
$S' \subseteq S$ is a finitely generated subgroup, then there
exists an ${\scriptstyle{{\bf H}}}\mathfrak{M}$-subgroup
$S'' \subseteq S$ with $S' \subseteq S''$. Since the
$kS''$-module $\mbox{res}_{S''}^GM$ is projective (in view
of the above discussion, regarding
${\scriptstyle{{\bf H}}}\mathfrak{M}$-subgroups of $G$),
this is also the case for its restriction
$\mbox{res}_{S'}^GM = \mbox{res}_{S'}^SN$ to $S'$; in
particular, the $kS'$-module $\mbox{res}_{S'}^SN$ is flat.
As $N$ is the filtered colimit of the $kS$-modules
$\mbox{ind}_{S'}^S\mbox{res}_{S'}^SN \in {\tt Flat}(kS)$,
where $S'$ runs through all finitely generated subgroups
of $S$, the $kS$-module $N$ is also flat. The projectivity
of $N$ follows now from the Benson-Goodearl theorem
\cite[Corollary 4.8]{BG}.

(ii) For any $\mathfrak{M}$-subgroup $F \subseteq S$, we
consider the $kF$-module $\mbox{res}_F^SN = \mbox{res}_F^GM$
and note that its restriction
$\mbox{res}_{H \cap F}^SN = \mbox{res}_{H \cap F}^GM =
 \mbox{res}_{H \cap F}^H\mbox{res}_H^GM$
to $H \cap F$ is projective (as a restriction of the
projective $kH$-module $\mbox{res}_H^GM$). We also note
that $(F,H \cap F,k)$ satisfies Moore's condition, whereas
the $\mathfrak{M}$-group $F$ satisfies Moore's conjecture;
it follows that $\mbox{res}_F^SN \in {\tt Proj}(kF)$. This
holds for any $\mathfrak{M}$-subgroup $F$ of the
$\Phi\mathfrak{M}$-group $S$ and hence
$\mbox{pd}_{kS}N < \infty$. Then, Proposition 5.2(i)
implies that $N \in {\tt Proj}(kS)$.

(iii) Working as in (ii) above, we conclude that the
$kF$-module $\mbox{res}_F^SN$ is projective and hence
flat for any $\mathfrak{M}$-subgroup $F$ of the
$\Phi_{flat}\mathfrak{M}$-group $S$. Therefore, it
follows that $\mbox{fd}_{kS}N < \infty$; an application
of Proposition 5.2(iv) shows that $N$ is actually projective.
\hfill $\Box$

\begin{Theorem}
Let $G$ be a group and $H \subseteq G$ a subgroup of finite
index, such that the triple $(G,H,k)$ satisfies Moore's
condition. If $M$ is a $kG$-module with
$\mbox{res}_H^GM \in {\tt Flat}(kH)$, then:

(i) $\mbox{res}_S^GM \in {\tt Flat}(kS)$ for any
${\scriptstyle{{\bf LH}}}\mathfrak{M}_{flat}$-subgroup
$S \subseteq G$,

(ii) $\mbox{res}_S^GM \in {\tt Flat}(kS)$ for any
$\Phi_{flat}\mathfrak{M}_{flat}$-subgroup $S \subseteq G$.
\end{Theorem}
\vspace{-0.05in}
\noindent
{\em Proof.}
We can follow the arguments in the proof of assertions
(i) and (ii) of Theorem 5.3, using Proposition 5.2(iii)
instead of Proposition 5.2(i). We note that there is no
need to appeal to the Benson-Goodearl theorem at the end
of the proof of (i). \hfill $\Box$

\begin{Theorem}
Let $G$ be a group and $H \subseteq G$ a subgroup of finite
index, such that the triple $(G,H,k)$ satisfies Moore's
condition. If $M$ is a $kG$-module with
$\mbox{res}_H^GM \in {\tt Inj}(kH)$, then:

(i) $\mbox{res}_S^GM \in {\tt Inj}(kS)$ for any
${\scriptstyle{{\bf LH}}}\mathfrak{M}_{inj}$-subgroup
$S \subseteq G$,

(ii) $\mbox{res}_S^GM \in {\tt Inj}(kS)$ for any
$\Phi_{inj}\mathfrak{M}_{inj}$-subgroup $S \subseteq G$.
\end{Theorem}
\vspace{-0.05in}
\noindent
{\em Proof.}
As we did in the proof of Theorem 5.3, for the subgroup
$S$ under consideration, we let $N = \mbox{res}_S^GM$
and note that
$\mbox{res}_{H \cap S}^SN = \mbox{res}_{H \cap S}^GM =
 \mbox{res}_{H \cap S}^H\mbox{res}_H^GM$.
Then, being a restriction of the injective $kH$-module
$\mbox{res}_H^GM$, the latter $k(H \cap S)$-module is
injective as well. We also note that $H \cap S$ has
finite index in $S$ and the triple $(S,H \cap S,k)$
satisfies Moore's condition.

(i) We first consider the case where $S$ is an
${\scriptstyle{{\bf H}}}\mathfrak{M}_{inj}$-subgroup of $G$
and proceed by induction on the ordinal $\alpha$, for which
$S \in {\scriptstyle{{\bf H}}}_{\alpha}\mathfrak{M}_{inj}$,
following the argument in the proof of \cite[Theorem 2.5]{BS1}.
If $\alpha =0$, then $S \in \mathfrak{M}_{inj}$ and hence $S$
satisfies the injective version of Moore's conjecture. The
$kS$-module $N$ is therefore injective (as its restriction
to $H \cap S$ is injective). For the inductive step, assume
that $\alpha >0$ and the result is true for all
${\scriptstyle{{\bf H}}}_{\beta}\mathfrak{M}_{inj}$-subgroups
of $G$ and all $\beta < \alpha$. As we noted in $\S $1.IV,
there exists a short exact sequence of $kS$-modules
\[ 0 \longrightarrow N \longrightarrow N^0
     \longrightarrow N^1 \longrightarrow \ldots
     \longrightarrow N^d \longrightarrow 0 , \]
where $d$ is the dimension of the $S$-CW-complex witnessing
that $S \in {\scriptstyle{{\bf H}}}_{\alpha}\mathfrak{M}_{inj}$
and each $N^i$ is a direct product of $kS$-modules of the form
$\mbox{coind}_F^S\mbox{res}_F^SN$, with $F$ an
${\scriptstyle{{\bf H}}}_{\beta}\mathfrak{M}_{inj}$-subgroup
of $S$ for some $\beta < \alpha$. Our inductive assumption
implies that the $kF$-module
$\mbox{res}_F^SN = \mbox{res}_F^GM$ is injective, so that
$\mbox{coind}_F^S\mbox{res}_F^SN \in {\tt Inj}(kS)$ for all
such subgroups $F$. Hence, $N^i$ is injective for all
$i=0,1, \ldots ,d$ and $\mbox{id}_{kS}N \leq d$. Finally,
an application of Proposition 5.2(ii) implies that the
$kS$-module $N$ is actually injective.

We now assume that $S$ is an
${\scriptstyle{{\bf LH}}}\mathfrak{M}_{inj}$-subgroup of $G$.
For any finitely generated subgroup $S' \subseteq S$, there
exists an ${\scriptstyle{{\bf H}}}\mathfrak{M}_{inj}$-subgroup
$S'' \subseteq S$ with $S' \subseteq S''$. Since the
$kS''$-module $\mbox{res}_{S''}^GM$ is injective (in view
of the above discussion, regarding
${\scriptstyle{{\bf H}}}\mathfrak{M}_{inj}$-subgroups of $G$),
this is also the case for its restriction
$\mbox{res}_{S'}^GM = \mbox{res}_{S'}^SN$ to $S'$. This holds
for any finitely generated subgroup $S' \subseteq S$ and hence
Lemma 5.6 below implies (by letting $H' = H \cap S$ therein)
that $N \in {\tt Inj}(kS)$.

(ii) We can follow the arguments in the proof of Theorem
5.3(ii), using Proposition 5.2(ii) instead of Proposition
5.2(i). \hfill $\Box$

\vspace{0.1in}

\noindent
The following Lemma was used in the proof of Theorem 5.5.
We have separated its proof from the main body of the proof
of that Theorem for expository reasons.

\begin{Lemma}
Let $S$ be a group, $H' \subseteq S$ a subgroup of finite
index and $N$ a $kS$-module. Assume that
$\mbox{res}_{H'}^SN \in {\tt Inj}(kH')$ and
$\mbox{res}_{S'}^SN \in {\tt Inj}(kS')$ for any finitely
generated subgroup $S' \subseteq S$. Then, $N \in {\tt Inj}(kS)$.
\end{Lemma}
\vspace{-0.05in}
\noindent
{\em Proof.}
First of all, we shall prove that
$\mbox{res}_{S'}^SN \in {\tt Inj}(kS')$ for any countable
subgroup $S' \subseteq S$. To that end, we fix a countable
subgroup $S' \subseteq S$ and express it as the union of an
increasing sequence $(S_i)_i$ of finitely generated subgroups.
If $M$ is a $kS'$-module, we consider the sequence of $kS'$-modules
$(M_i)_i$, where $M_i = \mbox{ind}_{S_i}^{S'}\mbox{res}_{S_i}^{S'}M$
for all $i$. The inclusions $S_i \hookrightarrow S_{i+1}$ induce on
$(M_i)_i$ the structure of a direct system, whose colimit is $M$.
Hence, there exists a short exact sequence of $kS'$-modules
\[ 0 \longrightarrow {\textstyle{\bigoplus_i}} M_i
     \longrightarrow {\textstyle{\bigoplus_i}} M_i
     \longrightarrow M \longrightarrow 0 . \]
Since
$\mbox{Ext}^n_{kS'} \! \left( M_i,\mbox{res}_{S'}^SN \right)
 \! = \mbox{Ext}^n_{kS_i} \! \left( \mbox{res}_{S_i}^{S'}M,
 \mbox{res}_{S_i}^SN \right) = 0$
for all $n>0$ and all $i$, the associated long exact sequence
of Ext-groups that corresponds to the above short exact sequence
of $kS'$-modules implies that
$\mbox{Ext}^n_{kS'} \left( M,\mbox{res}_{S'}^SN \right) = 0$
for all $n>1$. This is the case for any $kS'$-module $M$ and
hence $\mbox{id}_{kS'} \mbox{res}_{S'}^SN \leq 1$. Since
$H' \cap S'$ is a subgroup of finite index in $S'$ and the
restriction
$\mbox{res}_{H' \cap S'}^{S'}\mbox{res}_{S'}^SN =
 \mbox{res}_{H' \cap S'}^SN =
 \mbox{res}_{H' \cap S'}^{H'}\mbox{res}_{H'}^SN$
is injective, as a restriction of the injective $kH'$-module
$\mbox{res}_{H'}^SN$, Proposition 5.2(ii) implies that
$\mbox{res}_{S'}^SN \in {\tt Inj}(kS')$.

We shall now prove that $N$ is an injective $kS$-module,
using induction on the cardinality $\kappa$ of $S$. If
$\kappa \leq \aleph_0$, then the result follows from the
case considered above. If $\kappa$ is uncountable, then
we may express $S$ as a continuous ascending union of subgroups
$(S_{\lambda})_{\lambda < \kappa}$, each one having cardinality
$< \kappa$. Our induction hypothesis implies that the
$kS_{\lambda}$-module $\mbox{res}^S_{S_{\lambda}}N$ is injective
for all $\lambda$. (Indeed, $\mbox{res}^S_{S_{\lambda}}N$ is
injective when restricted to any finitely generated subgroup
of $S_{\lambda}$, as well as to the subgroup $H' \cap S_{\lambda}$,
which has finite index in $S_{\lambda}$.) Let $M$ be a $kS$-module
and define
$M_{\lambda} =
 \mbox{ind}_{S_{\lambda}}^S \mbox{res}^S_{S_{\lambda}}M$
for all $\lambda$. We note that
\[ \mbox{Ext}^n_{kS}(M_{\lambda},N) =
   \mbox{Ext}^n_{kS_{\lambda}} \! \left(
   \mbox{res}^S_{S_{\lambda}}M,\mbox{res}_{S_{\lambda}}^SN
   \right) = 0 \]
for all $\lambda$ and $n>0$. The ascending family of subgroups
$(S_{\lambda})_{\lambda < \kappa}$ induces a continuous direct
system of $kS$-modules $(M_{\lambda})_{\lambda < \kappa}$ with
surjective structure maps, whose colimit is $M$. The short exact
sequence of $kS$-modules
\[ 0 \longrightarrow K_{\lambda} \longrightarrow M_0
     \longrightarrow M_{\lambda} \longrightarrow 0 , \]
where $K_{\lambda}$ is the kernel of the structure map
$M_0 \longrightarrow M_{\lambda}$, is the $\lambda$-th term
of a continuous direct system of short exact sequences. The
colimit of the latter direct system of short exact sequences
is the short exact sequence of $kS$-modules
\begin{equation}
 0 \longrightarrow K \longrightarrow M_0
   \longrightarrow M \longrightarrow 0 .
\end{equation}
Here, $K$ is equal to the continuous ascending union of its
submodules $(K_{\lambda})_{\lambda < \kappa}$. Since
$K_{\lambda +1}/K_{\lambda}$ can be identified with
the kernel of the (surjective) structure map
$M_{\lambda} \longrightarrow M_{\lambda +1}$, we conclude that
$\mbox{Ext}^n_{kS}(K_{\lambda +1}/K_{\lambda},N)=0$ for all
$\lambda$ and $n>0$. Then, Eklof's lemma \cite[Theorem 7.3.4]{EJ2}
implies that $\mbox{Ext}^n_{kS}(K,N) = 0$ for all $n>0$; hence,
the short exact sequence (4) above implies that
$\mbox{Ext}^n_{kS}(M,N)=0$ for all $n>1$. Since this is the
case for any $kS$-module $M$, we conclude that
$\mbox{id}_{kS}N \leq 1$. Using once more Proposition 5.2(ii),
we conclude that $N$ is injective. \hfill $\Box$

\begin{Corollary}
(i) The class $\mathfrak{M}$ is closed under the
operations ${\scriptstyle{{\bf LH}}}$, $\Phi$ and
$\Phi_{flat}$.

(ii) The class $\mathfrak{M}_{flat}$ is closed under the
operations ${\scriptstyle{{\bf LH}}}$ and $\Phi_{flat}$.

(iii) The class $\mathfrak{M}_{inj}$ is closed under the
operations ${\scriptstyle{{\bf LH}}}$ and $\Phi_{inj}$.
\end{Corollary}
\vspace{-0.05in}
\noindent
{\em Proof.}
Theorem 5.3 implies that any group contained in
${\scriptstyle{{\bf LH}}}\mathfrak{M} \cup
 \Phi\mathfrak{M} \cup \Phi_{flat}\mathfrak{M}$
satisfies Moore's conjecture. In other words, we
have
${\scriptstyle{{\bf LH}}}\mathfrak{M} \cup
 \Phi\mathfrak{M} \cup \Phi_{flat}\mathfrak{M}
 = \mathfrak{M}$,
whence (i) holds. Analogously, Theorem 5.4 implies assertion
(ii) and Theorem 5.5 implies assertion (iii). \hfill $\Box$

\vspace{0.1in}

\noindent
As a consequence of Chouinard's work \cite{Chou}, all
finite groups satisfy Moore's conjecture. This result was
extended to the class of
${\scriptstyle{{\bf H}}}_1\mathfrak{F}$-groups in
\cite{ACGK}. In fact, all
${\scriptstyle{{\bf LH}}}\mathfrak{F}$-groups satisfy
Moore's conjecture; see \cite{AM} for the original
(projective) version of the conjecture and \cite{BS1}
for the injective version of the conjecture (which
implies the flat version, in view of
\cite[Proposition 2.6]{BS1}). Let
$\overline{\mathfrak{F}}$ be the smallest class of groups
which contains the class $\mathfrak{F}$ of finite groups
and is ${\scriptstyle{{\bf LH}}}$-closed, $\Phi$-closed
and $\Phi_{flat}$-closed; $\overline{\mathfrak{F}}$-groups
can be described hierarchically, as we explained in $\S 1$.IV.
The following Corollary is precisely Theorem B in the
Introduction. It is itself an immediate consequence of
Corollary 5.7(i), in view of the inclusion
$\mathfrak{F} \subseteq \mathfrak{M}$.

\begin{Corollary}
All $\overline{\mathfrak{F}}$-groups satisfy Moore's
conjecture over $k$, i.e.\
$\overline{\mathfrak{F}} \subseteq \mathfrak{M}$.
\end{Corollary}

\noindent
{\bf Remark 5.9.}
We may use Corollary 5.7(ii),(iii) and obtain analogues of
Corollary  5.8 for the flat and the injective versions of
Moore's conjecture.
\addtocounter{Lemma}{1}

\vspace{0.1in}

\noindent
{\em Acknowledgments.}
It is a pleasure to thank Sergio Estrada, Gregory Kendall,
Peter Kropholler and Olympia Talelli for useful comments
and suggestions on an earlier draft of this paper. Wei Ren
was supported by the National Natural Science Foundation of
China (No.\ 11871125).

\medskip

{\footnotesize \noindent Ioannis Emmanouil\\
Department of Mathematics, University of Athens,
Athens 15784, Greece \\
E-mail: {\tt emmanoui$\symbol{64}$math.uoa.gr}}

\medskip

{\footnotesize \noindent Wei Ren\\
School of Mathematical Sciences, Chongqing Normal
University, Chongqing 401331, PR China\\
E-mail: {\tt wren$\symbol{64}$cqnu.edu.cn}}

\end{document}